\documentclass[%
 aip,
 amsmath,amssymb,
preprint,%
]{revtex4-1}

\usepackage{graphicx}
\usepackage{dcolumn}
\usepackage{bm}
\usepackage[utf8]{inputenc}
\usepackage[T1]{fontenc}
\usepackage{mathptmx}
\usepackage{etoolbox}

\usepackage{booktabs}

\graphicspath{ {./figures/} }
\usepackage{subcaption}
\usepackage{hyperref}
\hypersetup{hidelinks}
\DeclareMathOperator*{\argmin}{arg\,min}

\makeatletter
\def\@email#1#2{%
 \endgroup
 \patchcmd{\titleblock@produce}
  {\frontmatter@RRAPformat}
  {\frontmatter@RRAPformat{\produce@RRAP{*#1\href{mailto:#2}{#2}}}\frontmatter@RRAPformat}
  {}{}
}%
\makeatother
\begin{document}

\preprint{ }

\title[ ]{Data-Driven Modeling and Forecasting of Chaotic Dynamics on Inertial Manifolds \\ Constructed as Spectral Submanifolds}
\author{Aihui Liu}
 \affiliation{ 
Institute for Mechanical Systems, ETH Zürich, Leonhardstrasse 21, 8092 Zürich, Switzerland.
}%

\author{Joar Ax{\aa}s}%
 \affiliation{ 
Institute for Mechanical Systems, ETH Zürich, Leonhardstrasse 21, 8092 Zürich, Switzerland.
}%

\author{George Haller}
  \email{georgehaller@ethz.ch}
\affiliation{ 
Institute for Mechanical Systems, ETH Zürich, Leonhardstrasse 21, 8092 Zürich, Switzerland.
}%

\date{\today}

\begin{abstract}
We present a data-driven and interpretable approach for reducing the dimensionality of chaotic systems using spectral submanifolds (SSMs). Emanating from fixed points or periodic orbits, these SSMs are low-dimensional inertial manifolds containing the chaotic attractor of the underlying high-dimensional system. The reduced dynamics on the SSMs turn out to predict chaotic dynamics accurately over a few Lyapunov times and also reproduce long-term statistical features, such as the largest Lyapunov exponents and probability distributions, of the chaotic attractor. We illustrate this methodology on numerical data sets including delay-embedded Lorenz and Rössler attractors, a nine-dimensional Lorenz model, a periodically forced Duffing oscillator chain, and the Kuramoto–Sivashinsky equation. We also demonstrate the predictive power of our approach by constructing an SSM-reduced model from unforced trajectories of a buckling beam, and then predicting its periodically forced chaotic response without using data from the forced beam.
\end{abstract}

\maketitle

\begin{quotation}
In high-dimensional systems, inertial manifolds containing the low-dimensional chaotic attractor play a significant role in understanding and predicting chaotic behaviors. Such low-dimensional, attracting inertial manifolds contain the essential dynamics of the system and therefore offer a tool for model reduction. We construct such inertial manifolds as spectral submanifolds (SSMs) emanating from an unstable steady state and develop SSM-reduced models on them. SSM-reduced models of chaotic systems give accurate short-term predictions and simultaneously reproduce long-term statistical properties. We further demonstrate the predictability of such SSM-reduced models on several numerical data sets including the delay-embedded Lorenz and Rössler attractors, forced Duffing oscillator chain, the Kuramoto–Sivashinsky equation and a periodically forced buckling beam. 
\end{quotation}

\section{Introduction}

Chaos is ubiquitous in a variety of systems, varying from the double pendulum to turbulent flows \cite{strogatz1994nonlinear,baker_gollub_1996,intro2DS}. Most of these chaotic systems are high-dimensional but their chaotic dynamics take place on low-dimensional attractors. The dynamics within these strange attractors provides an ideal reduced model for the full system, as long as appropriate localized coordinates can be constructed along the attractor. This typically requires the construction of a smooth, low-dimensional manifold into which the attractor is properly embedded. 

Based on Whitney's \cite{whitney} and Takens's \cite{takens} embedding theorems, several algorithms have been developed for attractor reconstruction and for learning appropriate coordinates near those attractors from delay-embedded data \cite{reConAttractor,attractorrecon,doi:10.1063/1.5039508}. In general, a higher embedding dimension preserves more information at the cost of increasing system dimensionality. To infer the optimal embedding dimension and delay embedding time lag, using average mutual information (AMI) \cite{ami} and the false nearest neighbor (FNN) method \cite{FNN} has become a standard practice. For high-dimensional time series, the dimensionality of the measurements is typically higher than that of the intrinsic system dynamics, requiring dimensional reduction for system identification and modeling. 


For this purpose, it is customary to first perform model reduction, then fit a low-dimensional chaotic model in the reduced space. Frequently used linear model reduction techniques, such as the principal component analysis (PCA) \cite{pca} or the dynamic mode decomposition (DMD) and its variants \cite{reviewDMD}, cannot capture intrinsically non-linearizable phenomena such as chaos. To counter this limitation, several nonlinear model reduction techniques have been developed, including kernel PCA and local linear embedding \cite{rom2,rom1}. Other popular methods for data-driven model reduction include autoencoders, deep learning and neural networks. These methods aim to extract hidden information with multiple levels of representation from nonlinear and complex dynamical systems. Various algorithms for black-box models have been developed, compared, and applied in various fields, including fluid dynamics and neural science \cite{Sangiorgio2021,Shahi2022}.

To improve the physical interpretability of models produced by these approaches, previous studies have utilized manifold learning \cite{Lee2020} or the concept of inertial manifolds \cite{inertial}. Inertial manifolds are smooth, finite-dimensional invariant manifolds that contain the global attractor and attract all solutions exponentially. In principle, therefore, inertial manifolds are the perfect tool for reducing complex systems. To utilize this tool, a recent method has been developed involving training neural networks to learn an atlas of low-dimensional inertial manifolds in high-dimensional spatio-temporal chaotic systems \cite{uwmadison2,uwmadison1}, which show significant promise.


Subsequent to model reduction, modeling chaos in the reduced coordinates is a further challenge. Such reduced models can be classified into local models and global models \cite{Abarbanel1996, a2}. Local models focus on capturing the dynamics within a small region of the state space for short-term prediction. An example of local chaos modeling is the $k$-th nearest neighbor (kNN) method \cite{knnFarmer}, which only relies on neighboring trajectory information to make forecasts. Although this method is inherently discrete and contains a large number of parameters, it has been proven to be practically useful when a significant amount of data is available \cite{Abarbanel1996}. In contrast, global models aim to capture the overall behavior across the whole phase space. Popular global modeling techniques include neural ODEs \cite{10.1063/5.0069536} and recurrent neural networks, such as LSTM \cite{LSTM1,DNN}. While these neural network models can make forecasts over short Lyapunov times and reproduce statistical properties of the system, they often require a large time domain for training data, extensive parameter tuning, and trial-and-error model selection.

A data-driven modeling alternative that offers simplicity and interpretability is the Sparse Identification of Nonlinear Dynamical Systems (SINDy) approach, which applies sparse regression to fit the right-hand side of a system of ODEs to the time-history of select observables that are assumed to form a closed, reduced dynamical system \cite{sindy}. While this methodology has been successful in several settings \cite{sindy2,Kaptanoglu2023}, its results are dependent on the choice of the selected function library, require an a priori knowledge of a reduced observable set, and suffer from sensitivity to noise. For increased interpretability, recent developments combine DMD or SINDy with physical principles \cite{piDMD,Silva2019DiscoveryOP}.

In this work, we introduce a data-driven model reduction method for chaotic dynamics using spectral submanifolds. A spectral submanifold (SSM) is the smoothest nonlinear continuation of a non-resonant spectral subspace of the linear part of a system at a stationary state \cite{SSM}. This stationary state can be a fixed point, periodic orbit, or quasiperiodic orbit. SSM-based model reduction has so far been used in the data-driven modeling of non-linearizable and deterministic dynamics in mechanical systems \cite{SSMLearn_appli}. SSM-reduction has also been successfully combined with model-predictive control of soft robots \cite{ssmControl1,ssmControl2} and augmented with a multivariate delay-embedding technique \cite{J_SSM}. In recent work, the mathematical theory behind SSM-reduction has been extended to fractional and mixed-mode SSMs \cite{mixedmodeSSM}, with the latter type of SSM also containing unstable modes in addition to the slowest stable modes.

This paper aims to further extend data-driven SSM-reduction to chaotic systems. Given an observable space that has a dimension higher than the Lyapunov dimension of the strange attractor, we seek to reconstruct a mixed-mode SSM that contains the chaotic attractor. Such an SSM then constitutes an explicitly constructed, very low-dimensional inertial manifold that is anchored to a stationary state. We illustrate on examples that mixed-mode, SSM-reduced models of chaotic systems give accurate short-term predictions while they simultaneously reproduce long-term statistical properties, such as the Lyapunov exponent and probability distribution density.

\section{Chaotic attractors on spectral submanifolds} 



\subsection{Spectral submanifold theory} \label{sec:SSMtheory}

We consider a dynamical system of the form
\begin{equation} \label{system}
    \dot{\bm{x}} = \bm{F} (\bm{x}) = \mathcal{A} \bm{x} + \bm{f} (\bm{x}) , \;\; \bm{x} \in \mathbb{R}^{n}, \;\; \mathcal{A} \in \mathbb{R}^{n \times n}, \;\; \bm{f} = \mathcal{O} (|\bm{x}|^{2}) \in \mathcal{C}^{\infty} .\;
\end{equation}
We assume that $\mathcal{A} \in \mathbb{R}^{n \times n}$ is diagonalizable and $\bm{x}=\bm{0}$ is a hyperbolic fixed point, i.e.,
\begin{equation}
    0 \notin \text{Re} \left[\text{spect}(\mathcal{A}) \right] ,
\end{equation}
where $\text{spect}(\mathcal{A}) = \{ \lambda _ 1, \lambda _ 2, ..., \lambda _{n} \}$ denotes the spectrum of the linear part of the dynamics at the origin. We further assume the following non-resonance conditions within $\text{spect}(\mathcal{A})$:
\begin{equation} \label{nonResonance}
    \lambda _j \neq \sum^{n}_{k=1} m_k \lambda _k , \;\; m_k \in \mathbb{N}, \;\; \sum^{n}_{k=1} m_k \geq 2, \;\; j = 1,...,n.
\end{equation}

We then select a $d$-dimensional spectral subspace $E \subset \mathbb{R}^{n}$ that is the direct sum of a family of real eigenspaces of $\mathcal{A}$. Without loss of generality, we assume that the spectrum of $\mathcal{A}$ within $E$, $\text{spect}(\mathcal{A}| _{E})$, contains $p$ purely real eigenvalues and $q$ pairs of complex conjugate eigenvalues, where $p+2q=d$. After a linear change of coordinates, the linear part of the system (\ref{system}) within the spectral subspace $E$ can be assumed diagonal. We denote the coordinates along the real eigenvectors by $\bm{u} \in \mathbb{R}^{p}$ and their corresponding Jordan block by 
\begin{equation} \label{eqa_p}
    \bm{A} = \text{diag} [ \lambda _1 ,..., \lambda _p] \in \mathbb{R}^{p \times p} .
\end{equation}
Similarly, we denote the coordinates along the complex conjugate $q$ pairs of eigenvectors by $\bm{z} \in \mathbb{C}^{q}$, and their real Jordan block by
\begin{equation} \label{eqa_q}
    \bm{B}_{i} = \begin{pmatrix}
                 \alpha _i & \omega _i \\
                 -\omega _i & \alpha _i  \end{pmatrix}
                 \in \mathbb{R}^{2 \times 2}, \;\; i = 1,..., q .
\end{equation}

We finally assume that the rest of the spectrum of $\mathcal{A}$ outside $\text{spect}(\mathcal{A}| _{E})$ comprises $r$ purely real eigenvalues and $s$ pairs of complex conjugate eigenvalues that correspond to the coordinates 
\begin{equation}
    \bm{v} \in \mathbb{R}^{r}, \;\; \bm{w} \in \mathbb{C}^{s},
\end{equation}
respectively, within our final set of coordinates $(\bm{u},\bm{z}, \bm{v}, \bm{w})$.

Spectral submanifolds (SSMs) are invariant manifolds $\mathcal{W}(E)$ that serve as nonlinear continuations of a $(p+2q)$-dimensional spectral subspace
\begin{equation}
    E = \{(\bm{u},\bm{z}, \bm{v}, \bm{w}) : \bm{v} = \bm{0}, \bm{w} = \bm{0} \}
\end{equation}
of the linearized dynamics at a fixed point. When the linear part of system (\ref{system}) is asymptotically stable and the non-resonance conditions (\ref{nonResonance}) hold, then there is a unique, primary $\mathcal{W}(E)$ that is as smooth as the full system \cite{SSM,SSMLearn}. Recent work by \cite{mixedmodeSSM} extended this result to general hyperbolic fixed points with general (like-mode or mixed-mode) spectral subspaces. (Like-mode spectral subspaces contain only eigenspaces of the same stability type, in contrast to their mixed-mode counterparts.) In this paper, we focus only on the unique, $\mathcal{C}^{\infty}$-smooth primary SSMs, $\mathcal{W}(E)$, that are tangent to mixed-mode spectral subspaces $E$ of hyperbolic fixed points.

By the existence theory of such primary SSMs \cite{mixedmodeSSM}, for any positive integer $\mathcal{K} \geq 2$, there exists a unique set of coefficients stored in a matrix $\bm{G} $ such that in the $(\bm{u},\bm{z},\bm{v},\bm{w})$ coordinates introduced above, the SSM, $\mathcal{W}(E)$, tangent to a spectral subspace $E$ can locally be written near $\bm{x} = \bm{0}$ as a $\mathcal{C}^{\infty}$ graph over the $(\bm{u},\bm{z})$ variables:
\begin{equation}
    \begin{aligned}
        (\bm{v},\bm{w})^{\mathrm{T}}
        & = \mathcal{G}(\bm{u}, \bm{z}) = \bm{G} \cdot (\bm{u}, \bm{z}, \overline{\bm{z}}) ^{1 : \mathcal{K}} + \mathcal{O} \left( |(\bm{u},\bm{z})|^{\mathcal{K}} \right) , \\
        \bm{G} & = [\bm{G}_1 , \bm{G}_2, \dots, \bm{G}_{\mathcal{K}}] \in \mathbb{C}^{(n-d) \times \sum^{\mathcal{K}}_{k=1} d_k }, \;\; \bm{G}_k \in \mathbb{C}^{(n-d) \times d_k}. \;
    \end{aligned}
\end{equation}

Here the superscript $(\cdot)^{l:r}$ denotes a vector of all monomials of the variables listed in $(\cdot)$ from order $l$ to order $r$ and $d_k$ is the total number of $d$-variate monomials composed of those variables at order $k$. For instance, if $p = 2, \;q = 0 $ in Eqs. (\ref{eqa_p}) and (\ref{eqa_q}), and the polynomial order we select to approximate the SSM is $\mathcal{K} = 3$, then we have $(\bm{u}, \bm{z}, \overline{\bm{z}}) = [u_1, u_2]^{\mathrm{T}}$, $d_2 = 3$ and $d_3 = 4$, generating the set of monomials
\begin{equation}
    (\bm{u}, \bm{z}, \overline{\bm{z}})^{2:3} = [u_1^2, u_1 u_2, u_2^2, u_1^3, u_1^2 u_2, u_1 u_2^2, u_2^3]^{\mathrm{T}}, 
\end{equation}
and the coefficient matrix $\bm{G}_{2:3} \in \mathbb{R}^{ (n-2) \times (3+4)}$ is real. 

The reduced dynamics on the SSM, $\mathcal{W}(E)$, are given by
\begin{equation}
    \begin{pmatrix} \dot{\bm{u}} \\ \dot{\bm{z}} \end{pmatrix} = 
    \begin{pmatrix}
        \bm{A} & 0 & 0 & 0 \\ 0 & \bm{B}_1 & 0 & 0 \\ 0 & 0 & \ddots & 0 \\ 0 & 0 & 0 & \bm{B}_q
    \end{pmatrix}
    \begin{pmatrix} \bm{u} \\ \bm{z}     \end{pmatrix} + 
    f^{\bm{u},\bm{z}} (\bm{u}, \bm{z}, \mathcal{G}(\bm{u}, \bm{z})),
\end{equation}
where $f^{\bm{u},\bm{z}}$ denotes the $(\bm{u}, \bm{z})$ coordinate component of $f = ( f^{\bm{u},\bm{z}}, f^{\bm{v}}, f^{\bm{w}} )$.

The open source \textsc{Matlab} software \texttt{SSMTool} \cite{tool} can compute a Taylor series approximation of the SSM, $\mathcal{W}(E)$, up to any required order when system (\ref{system}) is explicitly known. When system (\ref{system}) is only known from data, then the open source \texttt{SSMLearn} \cite{SSMLearn} algorithm can approximate $\mathcal{W}(E)$ purely from observed trajectories.

\subsection{SSM reconstruction via delay embedding} \label{sec:delay}

In practice, physical systems in which all state variables can be simultaneously monitored are rare. This necessitates the use of invariant manifold reconstruction methods relying only on a limited number of observables. Specifically, using the Takens delay-embedding theorem \cite{takens,embedo}, it is possible to embed an invariant manifold of dimension $d$ in the space of $2d+1$ (or higher) time-shifted samples of a single, generic scalar observable. Recently, Ref. \cite{J_SSM} discussed the use of delay embedding to recover specifically SSMs and their tangent spaces in such observable spaces.

To recall this methodology, we let $\mathcal{W}$ be a $d$-dimensional spectral submanifold of (\ref{system}), and denote the reduced dynamics on $\mathcal{W}$ by $\dot{\bm{\eta}} = \bm{R}(\bm{\eta})$, $\bm{\eta} \in \mathbb{R}^{d}$, which generates the reduced flow map $\bm{R}^{t} : \bm{\eta_0} \mapsto \bm{\eta}(t,\bm{\eta_0})$. Then there exists a smooth SSM parametrization $\bm{M}: \mathbb{R}^{d} \mapsto \mathbb{R}^{n}$, such that the reduced flow map $\bm{R}^{t}$ is conjugate to the restriction of the full flow map $\bm{F}^{t}: \mathbb{R}^n \mapsto \mathbb{R}^n$ of system (\ref{system}) to $\mathcal{W} \subset \mathbb{R}^n$:
\begin{equation} \label{coj1}
    \bm{F}^{t} \circ \bm{M} = \bm{M} \circ \bm{R}^{t}  .
\end{equation}

We then consider time-resolved, scalar-valued observations $s (t) = \zeta (\bm{x}(t))$ of system (\ref{system}) from a generic observable $\zeta: \mathbb{R}^{n} \mapsto \mathbb{R} $. We define a delay coordinate map with $m$ delays: $\bm{\Psi}  : \mathbb{R}^{n} \mapsto \mathbb{R}^{m}$ by stacking $m$ subsequent entries of the time series $s (t)$, sampled uniformly in time with sampling time $T_s$, in a vector:
\begin{equation}
    \bm{y}(t) = \bm{\Psi} (\bm{x}(t)) = [ s (t), s (t+T_s), \cdots , s ( t+(m-1) T_s ) ]^{\mathrm{T}} \in \mathbb{R}^{m} . 
\end{equation}

The map $\bm{\Psi} $ is thus a mapping from the phase space of system (\ref{system}) to the delay embedded space $\mathbb{R}^{m}$. Let $\bm{F}_{\Psi} : \mathbb{R}^{m} \mapsto \mathbb{R}^{m}$ be the flow map induced by the flow of system (\ref{system}) on this delay embedding space. By construction, if the observation map $\zeta$ maps $\bm{0}$ to $0$, the $\bm{x} = \bm{0}$ equilibrium point of system (\ref{system}) is mapped into a fixed point $\bm{y} = \bm{0}$, i.e., $\bm{F}_{\Psi} (\bm{0}) = \bm{0}$. According to the Takens embedding theorem \cite{takens,embedo}, for a generic observable function $\zeta$ and under certain nondegeneracy conditions on the sampling time $T_s$, if $m>2d$, then the delay coordinate map $\bm{\Psi}$ restricted to $\mathcal{W}$ provides an embedding of $\mathcal{W}$ into $\mathbb{R}^{m}$. As a consequence, $\Tilde{\mathcal{W}} = \bm{\Psi}(\mathcal{W})$ is diffeomorphic to the manifold $\mathcal{W}$. 

The delay-embedded dynamics $\bm{F}_{\Psi}$ on $\Tilde{\mathcal{W}}$ and the reduced dynamics restricted to the invariant manifold $\mathcal{W}$ are conjugate, i.e.,
\begin{equation} \label{coj2}
     \bm{F}_{\Psi} \circ \bm{\Psi} =  \bm{\Psi} \circ \bm{F}^{t} | _{\mathcal{W}}  , 
\end{equation}
which provides the theoretical basis for identifying an SSM attached to the $\bm{y} = \bm{0}$ origin of the $m$-dimensional delay-embedding space. 

For dynamical systems with (generally non-smooth) global attractors, inertial manifolds have been defined as finite-dimensional, forward-invariant manifolds that are at least Lipschitz and contain the global attractor \cite{inertial,Temam1990}. Therefore, attracting SSMs with an attractor in their reduced dynamics function as inertial manifolds, at least locally near a stationary state. If the attractor within the SSM is globally unique, then an SSM is, in fact, an inertial manifold. Unlike general inertial manifolds, however, SSMs can be precisely located and computed up to arbitrary order of accuracy, which makes them the perfect tools for a mathematically precise reduction of high-dimensional systems to very low-dimensional models. 

In order to contain the attractor, the dimension of the SSM should be larger than the Hausdorff dimension of the attractor \cite{robinson_2010}. In a data-driven setting with no prior information on the attractor dimension, the necessary embedding dimension is often estimated via the false nearest neighbors (FNN) method \cite{strogatz1994nonlinear,kantz_schreiber_2003}. This method is based on the idea that if the embedding dimension is too low, neighboring points in the embedding space may appear closer than they actually are, leading to false neighbors. Specifically, the FNN algorithm examines data points in the delay-embedding space, and checks if they are indeed close to each other in all directions. The percentage of the false nearest neighbors goes down when the embedding dimension is increased. When the percentage becomes lower than a predefined value, the data disentangle in the embedding space, thus giving a correct approximation for the attractor dimension. The SSMs we construct will have a dimension of $d$ that is at least as large as the estimate obtained from the FNN algorithm.

\section{Learning spectral submanifolds from data}

According to the SSM theory discussed in Section \ref{sec:SSMtheory}, the primary SSM, $\mathcal{W}(E)$, tangent to the spectral subspace $E$ at the origin, can locally be written as a polynomial expansion over $E$. The open source \textsc{Matlab} and Python implementations of \texttt{SSMLearn} \cite{ssmlearn_code} compute the SSM parametrization and its reduced dynamics purely from observed data \cite{SSMLearn,SSMLearn_appli}. In this section, we discuss the computation of the SSM-reduced dynamics and predictions from such reduced models of the full system dynamics.

\subsection{Data-driven SSM-based dimensionality reduction}

To provide a clearer explanation, in this section, we only consider $q=0$, i.e., the chosen spectral subspace $E$ only contains $p$ real eigenvalues. Our procedure, however, extends directly to spectral subspaces with complex eigenvalues, such as those arising in \cite{mixedmodeSSM}. The applications of numerical datasets in Section \ref{sec:application} contain examples with general choices of spectral subspaces, i.e., with both $q=0$ and $q \neq 0$. 


Let the coordinate $\bm{y} \in \mathbb{R}^{\rho}$ denote data obtained from system (\ref{system}) either via direct observations or delay-embedding (see Section \ref{sec:delay}). Also, let $\bm{\eta} \in \mathbb{R}^{d}$ denote coordinates along the corresponding spectral subspace $E$ of $D\bm{\Psi}(E)$, respectively. As a graph over this spectral subspace, the embedded SSM can be approximated by a polynomial expansion. According to the $\mathcal{C}^{\infty}$ version of the SSM theory outlined in Section \ref{sec:SSMtheory}, for any positive integer order $\mathcal{K} \geq 2$, the coefficients of this polynomial expansion are unique. We, therefore, seek a parametrization $\bm{M} : \mathbb{R}^{d} \mapsto \mathbb{R}^{\rho}$ of the embedded SSM, $\Tilde{\mathcal{W}}(E) \subset \mathbb{R}^{\rho}$, as a $\mathcal{K}$-th order polynomial in the form
\begin{equation} \label{Mpara}
\begin{aligned}
        \bm{M} ( \bm{\eta} ) &  = \bm{V} \bm{\eta}^{1:\mathcal{K}} = \bm{V}_1 \bm{\eta} + \bm{V}_{2:\mathcal{K}} \bm{\eta}^{2:\mathcal{K}} , \\
        \bm{V} & = [\bm{V}_1 , \bm{V}_2, \dots, \bm{V}_{\mathcal{K}}] \in \mathbb{R}^{\rho \times  \sum^{\mathcal{K}}_{k=1} d_k }, \;\; \bm{V}_k \in \mathbb{R}^{\rho \times d_k},
\end{aligned}
\end{equation}
with the reduced coordinates $\bm{\eta} = \bm{V}_1^{\mathrm{T}} \bm{y}$, where the matrix $\bm{V}_1 \in \mathbb{R}^{\rho \times d}$ has orthonormal columns that span the tangent space of the manifold $\Tilde{\mathcal{W}}(E) \subset \mathbb{R}^{\rho}$ at $\bm{y} = \bm{0}$. 

To find the constant in the above parametrization from data, we use \texttt{SSMLearn} to solve the constrained minimization problem
\begin{equation} \label{minimization}
\begin{aligned}
\bm{V}^{*} = [\bm{V}_1^{*} , \bm{V}_{2:\mathcal{K}}^{*}] = &  \argmin_{\bm{V}_1 , \bm{V}_{2:\mathcal{K}}} \quad  \sum _{j} ||  \bm{y}_j - \bm{V}_1 \bm{V}_1^{\mathrm{T}} \bm{y}_j - \bm{V}_{2:\mathcal{K}} \bm{\eta}^{2:\mathcal{K}} ||  , \\
& \bm{V}_1^{\mathrm{T}}\bm{V}_1 = \bm{I} , \quad \bm{V}_1^{\mathrm{T}} \bm{V}_{2:\mathcal{K}} = \bm{0} , \\
\end{aligned}
\end{equation}
where the last constraint represents a basic nonlinear extension of the principal component analysis.

To quantify the accuracy of this SSM approximation, we define the \textit{invariance error} as follows. From a specific observation $\bm{y}_{i} \neq \bm{0}$, we compute the corresponding reduced coordinate $\bm{\eta}_{i} = \bm{V}_1 ^{\mathrm{T}} \bm{y}_{i} $ as the projection of $\bm{y}_{i}$ onto the tangent space $\bm{V}_1$. By lifting $\bm{\eta}_{i}$ to the full embedding space via the manifold parametrization (\ref{Mpara}) and comparing it with the original data point $\bm{y}_{i}$, we obtain a measurement of the manifold fitting accuracy for a given choice of the approximation order $\mathcal{K}$. For a total of $P$ observations, we denote the lifted data points by $\bm{M}(\bm{\eta}) = \bm{M}( \bm{V}_1^{\mathrm{T}}\bm{y}_{i} )$ and define the normalized invariance error as
\begin{equation} \label{NMFE}
    \text{Invariance Error} = \frac{1}{||\underline{\bm{y}}||} \frac{1}{P} \sum_{i = 1}^{P} || \bm{y}_{i} - \bm{M}(\bm{V}_1^{\mathrm{T}}\bm{y}_{i}) || , \qquad || \underline{\bm{y}} || = \max_{i} || \bm{y}_i ||.
\end{equation}
A larger manifold parametrization order $\mathcal{K}$ will generally decrease the invariance error, but excessively large $\mathcal{K}$ values will lead to overfitting. In the examples from this paper, $\mathcal{K}$ is often chosen at the local minimum of the invariance error. The value of the invariance error depends on specific problems, but keeping it below $0.5\%$ suffices in our experience.

\subsection{Prediction from SSM-reduced order models} \label{sec:prediction}

We obtain a model from the reduced dynamics on the SSM by fitting the right-hand side of an ODE to the numerically determined, SSM-reduced vector field in the reduced coordinates. In this paper, we employ two alternative approximations for the reduced dynamics models: a polynomial fit and a nearest neighbor interpolation.  

In our first prediction method, we approximate the reduced dynamics in the form
\begin{equation} \label{rddynamics}
    \dot{\bm{\eta}} = \bm{R}_1 \bm{\eta} + \bm{R}_{2:\mathcal{K}} \bm{\eta}^{2:\mathcal{K}} ,
\end{equation}
where the coefficient matrix $ \bm{R}  = [\bm{R}_1, \bm{R}_{2:\mathcal{K}}]  \in \mathbb{R} ^{d \times \sum^{\mathcal{K}}_{k=1} d_k }$ is obtained via regression as
\begin{equation}
    \bm{R}^{*} =\argmin _{\bm{R} } \; \sum _{j} || \dot{\bm{\eta}}_j - \bm{R} \bm{\eta}_{j}^{1:\mathcal{K}} || .
\end{equation}

We then diagonalize the linear part of the reduced dynamics (\ref{rddynamics}) by a linear transformation $\bm{\xi} = \bm{W}^{-1} \bm{\eta}$, where $\bm{\xi} \in \mathbb{R}^{d}$ denotes modal coordinates and $\bm{W}\in \mathbb{R}^{d \times d}$ is the matrix of the eigenvectors of $\bm{R}_1$. In those new coordinates, the ODE (\ref{rddynamics}) becomes
\begin{equation} \label{modaldynamics}
    \dot{\bm{\xi}} = \bm{\Lambda} \bm{\xi} + \bm{N}_{2:\mathcal{K}} \bm{\xi}^{2:\mathcal{K}},
\end{equation}
where $\bm{\Lambda} \in \mathbb{R}^{d \times d}$ is a diagonal matrix containing the eigenvalues of $\bm{R}^{*}$, and $\bm{N}_{2:\mathcal{K}} \in \mathbb{R}^{d \times \sum^{\mathcal{K}}_{k=2} d_k }$ is the coefficient matrix of nonlinear monomials of $\bm{\xi}$ from order $2$ to $\mathcal{K}$. The coefficient matrix $\bm{N}_{2:\mathcal{K}}$ can be specified via a recursive sequence of normal form transformations that preserve the dynamics, as explained in \cite{SSMLearn}. Here we do not perform this step because subsequent normal form transformations are generally defined on increasingly smaller neighborhoods of $\bm{y} = \bm{0}$. This would limit our ability to capture larger-sized attractors on the SSM.

A global polynomial representation of the reduced dynamics may require very high-order polynomials, and yet fail to reproduce the chaotic dynamics \cite{Abarbanel1996}. Also, the lack of training data outside of the attractor may lead to limited data structures in monomial spaces, which makes it difficult to fit a globally stable polynomial model to vector fields with chaotic trajectories. As an alternative, for complicated strange attractors, we turn to the simplest and earliest method of local forecasting \cite{local_pr,10.1016/j.camwa.2009.10.019}: we fit maps locally to data points, and then use the resulting collection of local polynomials to form a global model. 

This idea brings us to our second prediction method, the local $k$-th nearest neighbor (kNN) approach. The kNN prediction method \cite{knnFarmer,kNN2} identifies the $k$ nearest neighbors in the training data set and then, by using various interpolation or regression methods, makes local predictions based on those neighbors \cite{knn3}. This method is suitable for predicting highly nonlinear or even chaotic behavior \cite{knn4}. In this paper, for a prediction based at a point $\bm{\eta}$, we first search for its $k$-th nearest neighbor $\bm{\eta}_j$, $j = 1,...,k$ in the training data set that minimizes the Euclidean norm $\| \bm{\eta} - \bm{\eta}_j \| _2$. We then use a zeroth-order linear fit to give the prediction $\bm{\eta}^{*}$ at the next time step, i.e., define the weight $w_j$ as
\begin{equation}
    w_j = \frac{ \| \bm{\eta}-\bm{\eta}_j \|_2 } { \sum^{k}_{j = 1} \| \bm{\eta} - \bm{\eta}_j \|_2 },
\end{equation}
and take the weighted sum
\begin{equation}
    \bm{\eta}^{*} = \sum^{k}_{j = 1}  w_j \| \bm{\eta} - \bm{\eta}_j \|_2 .
\end{equation}

Although the kNN method, like other local chaos prediction methods, has the limitation of discontinuities, large numbers of adjustable parameters, and the ability to predict only within the training data range, as long as the training data is sufficiently dense on the attractor, the method still has good predictive power. 

In this paper, we only use polynomial regression and the kNN method to fit the reduced dynamics. After obtaining the data sets in reduced coordinates, however, it is also possible to combine SSM-based model reduction with other forecasting techniques. We discuss this further in the Conclusions.



\subsection{Reconstruction criteria}

For non-chaotic dynamics, the prediction errors from SSM-reduced models are quantified by the normalized mean-trajectory-error (NMTE) \cite{SSMLearn}, which is the normalized difference between the test trajectory and the predicted trajectory. Chaotic dynamics, however, have sensitive dependence on initial conditions, which implies that the predictions for individual trajectories will invariably fail for longer times. Indeed, even a well-constructed reduced model can only be expected to capture statistical measures of the chaotic dynamics in long-term predictions. Accordingly, in this paper, we will evaluate the ability of SSM-reduced models to reproduce both the Lyapunov exponent and the probability density distribution of the attractors. 

First, we consider the Lyapunov characteristic exponents \cite{Loriginal,lyapunov}, the most often used indicators of chaotic flows. Let $\bm{x}(t;\bm{x}_0)$ denote a trajectory of system (\ref{system}) starting from the initial condition $\bm{x}_0$ at $t = 0$. An initial perturbation $\bm{\varepsilon} (0) = \bm{\varepsilon} _{0}$ to $\bm{x_0}$ will evolve along $\bm{x}(t;\bm{x}_0)$ as $\bm{\varepsilon}_{t} = \bm{\Phi }(t;\bm{x}_0) \bm{\varepsilon} _{0}$, where $ \bm{\Phi }(t;\bm{x}_0) $ is the fundamental matrix solution of the linearized system (\ref{system}) along the trajectory $\bm{x} (t;\bm{x}_0)$. The maximal exponential rate of expansion along $\bm{x} (t;\bm{x}_0)$ is then the maximal Lyapunov exponent
\begin{equation}
    \lambda _{\text{max}} (\bm{x}_0) = \lim _{t \rightarrow \infty} \frac{1}{2t} \log  \lambda _{\text{max}} \left[ \bm{\Phi}^{\mathrm{T}} (t;\bm{x}_0) \bm{\Phi} (t;\bm{x}_0) \right]
\end{equation}
if the limit exists.

As the maximum rate of separation of nearby trajectories, $\lambda _{\text{max}} (\bm{x}_0)$ defines a local measure of predictability for a dynamical system. The Lyapunov exponents are independent of the initial condition $\bm{x}_{0}$ and depend only on the trajectory $\bm{x}(t;\bm{x}_0)$. If $\bm{x}(t;\bm{x}_0)$ is contained in a chaotic (or strange) attractor, then $\lambda _{\text{max}}(\bm{x}_0)$ is constant for almost all trajectories along the attractor \cite{chaoticDynamics}.

One way to approximate $\lambda _{\text{max}} (\bm{x}_0)$ is to compute the exponent of the rate at which trajectories from nearby initial conditions diverge from the trajectory $\bm{x}(t;\bm{x}_0)$, i.e., calculate
\begin{equation}
    \lambda (t; \bm{x}_0 , \bm{\epsilon}_0) = \log \frac{ | \bm{x}(t; \bm{x}_0+\bm{\epsilon}_0) - \bm{x}(t;\bm{x}_0) | }{ |\bm{\epsilon}_0| }
\end{equation}
for a randomly chosen initial perturbation $\bm{\epsilon}_0$ to $\bm{x}_0$ and for large enough times $t$. For a chaotic system, $\lambda (t; \bm{x}_0, \bm{\epsilon}_0)$ is approximately a linear function of time \cite{Rosenstein1993APM}. When a good linear fit to $\lambda (t; \bm{x}_0, \bm{\epsilon}_0)$ is not possible, the system is generally not chaotic \cite{lyapunov2}. While $\lambda (t; \bm{x}_0, \bm{\epsilon}_0) /t$ depends on the initial condition $\bm{x}_0$ and the perturbation $\bm{\epsilon}_0$, this ratio will converge to $\lambda _{\text{max}} (\bm{x}_0)$ for generic choices of $\bm{x}_0$ and $\bm{\epsilon}_0$ \cite{chaoticDynamics}. Based on this observation, we randomly pick two nearby initial conditions on chaotic attractors, and use a linear fit to their temporal separation rate to estimate the Lyapunov exponent $\lambda _{\text{max}}$ of the attractor.

A related measure of predictability is the Lyapunov time, defined as the reciprocal of $\lambda _{\text{max}}$ of the attractor \cite{gaspard_1998}, i.e.,
\begin{equation}
    \text{Lyapunov}\; \text{Time} = \frac{1}{\lambda _{\text{max}}}.
\end{equation}
This quantity gives the characteristic time scale of chaotic dynamics. We call a prediction a short-time prediction if it is carried out over times that are low-order multiples of the Lyapunov time.  

Along with the Lyapunov exponent, the probability density \cite{ruelle_1989,10.1214/ss/1177011444} is also an invariant property of attractors and hence is often used for comparisons between the original and modeled chaotic attractors \cite{pdf1,PhysRevA.22.1198}. If a chaotic attractor contains the orbit $\bm{x}(t)$, its probability measure $\sigma$ is defined as the asymptotic average of Dirac deltas along $\bm{x}(t)$ \cite{RevModPhys.57.617}, i.e., 
\begin{equation}
    \sigma  = \lim _{T \to \infty} \frac{1}{T} \int _{0} ^{T} \delta _{x(t)} dt ,
\end{equation}
where the integral of $\delta _{x(t)}$ over the entire trajectory is equal to one.

On a chaotic attractor, the invariant measure $\sigma$ is independent of the initial conditions. To approximate the probability density of the chaotic attractors, we will take a sampling time that is long enough for $\sigma$ to converge. Along each axis of the system, we can either count the number of data points between intervals and plot the histogram, or estimate the probability distribution $\sigma$ based on a normal kernel function \cite{a2}. In our examples, we will compare the probability density of the original data points with those of the predictions made by SSM-reduced model using the same kernel function. The closeness of the two plots indicates how statistically accurate our long-term SSM-based predictions are.

In summary, we will test SSM-based prediction quality for short times via trajectory error calculations, and for long times via Lyapunov exponent and probability density distribution calculations.

\section{Applications to chaotic systems} \label{sec:application}

In this section, we perform SSM-based model reduction on four numerical data sets containing chaotic attractors.

In our first two examples, we observe only a scalar time series and delay-embed an SSM containing the attractor. In our subsequent examples, we assume that all phase space variables can be fully observed. For the first delay-embedded Lorenz attractor, we use a global polynomial representation for the reduced dynamics. For the other examples, we adopt the kNN nearest neighbor approach.

For autonomous ODE examples, the SSM-reduction outlined in Section \ref{sec:SSMtheory} will be directly applicable. We also explore one non-autonomous (time-periodic) ODE example in Section \ref{sec:duffing}, where we will carry out the discrete version of the same procedure for the Poincaré map of the system. Specifically, we perform model reduction on the Poincaré map, identifying an SSM emanating from its unstable fixed point.

We further include the example of the delay-embedded Rössler attractor in Section \ref{sec:rossler}, in which the SSM is attached to the unstable subspace corresponding to a complex conjugate pair of eigenvalues. The applications of SSM-based model reduction range from classic chaotic attractors to low-dimensional attractors in high-dimensional systems, such as Kuramoto-Sivashinski equation in Section \ref{sec:kse} and the finite element beam model in Section \ref{sec:vkbeam}.

\subsection{Delay-embedded Lorenz attractor} \label{sec:3Dlorenz}

As our first example, we consider the Lorenz model \cite{Lorenz} 
\begin{equation} \label{lorenz_eqa}
    \begin{aligned}
    \dot{x} &= \sigma (y-x) ,\\
    \dot{y} &= x(\rho -z) -y  ,\\
    \dot{z} &= xy-\beta z,  \;
    \end{aligned}
\end{equation}
with the classic parameter values $\sigma = 10$, $\rho = 28$ and $ \beta = 8/3 $. This system has an unstable fixed point at $(x,y,z) = \bm{0}$ and hence falls in the setting of this paper. We will only use a single scalar signal in our analysis and hence our methodology is fully unaware of the dimension and form of the ODE (\ref{lorenz_eqa}). As an SSM containing a chaotic attractor must be at least three-dimensional, Takens's theorem guarantees success (with probability one) in embedding the SSM in an at least seven-dimensional space. We then seek to reconstruct the Lorenz attractor within the 3D SSM.

To compare the SSM-reduced model results with the broadly used SINDy algorithm mentioned in the Introduction, we will adopt the same setting as in the \textsc{Matlab} code given by \cite{sindy}. Accordingly, we will observe only the $x$ coordinate of system (\ref{lorenz_eqa}) to generate a scalar observable time series. We generate two trajectories (one for training and the other for testing) over the interval $(0,100)$ with a time step of $0.001$. The test data starts with the classic initial condition $(-8,8,27)^{\mathrm{T}}$, while the initial condition of the training data is slightly separated $(-8,8,27+\epsilon)^{\mathrm{T}}$, where $\epsilon = 10^{-9}$. After truncating the data over the first short period of the time interval $(0,1)$, each of these two individual trajectories lies on the low-dimensional attractor of system (\ref{lorenz_eqa}).

We use FNN to estimate the dimension $d$ of the SSM that prevails in the data. The different choices of model reduction parameters and corresponding manifold fitting results are shown in Table \ref{3D_dim}. The high FNN percentage, $31.4\%$, in the first row of Table \ref{3D_dim} indicates that the SSM containing the strange attractor has a dimension higher than two. This prompts us to select the SSM dimension as $d= 3$, Based on the Takens embedding theorem, a 3D SSM can generically be delay-embedded in a $\rho = 2 \times d +1 = 7$ dimensional space. Table \ref{3D_dim} shows the accuracy of SSM reconstruction in this setting under different polynomial approximation orders for the 3D SSM.

\begin{table}[h]
\centering
\begin{tabular}{ccccc}
\toprule
SSM dimension ($d$)           & 2 & 3 & 3 & 3 \\ \midrule
FNN percentage           & 31.4\%  & 0.0\%  & 0.0\%  & 0.0\%  \\ \midrule
Delay embedding dimension  ($\rho$)    & 5 & 7 & 7 & 7 \\ \midrule
Manifold expansion order ($\mathcal{K}$) & - & 2 & 3 & 4 \\ \midrule
Invariance Error         & - & 3.85\%  & 1.25\%  & 1.45\%  \\ \bottomrule
\end{tabular}
\caption{Different choices of SSM dimension, false nearest neighbor percentage (FNN), order of manifold parametrization, and the corresponding invariance error. Due to the high FNN percentage for $d=2$, we did not compute an SSM approximation for that case.}
\label{3D_dim}
\end{table}

\begin{figure}[h]
    \centering
    \includegraphics[width=0.5\textwidth]{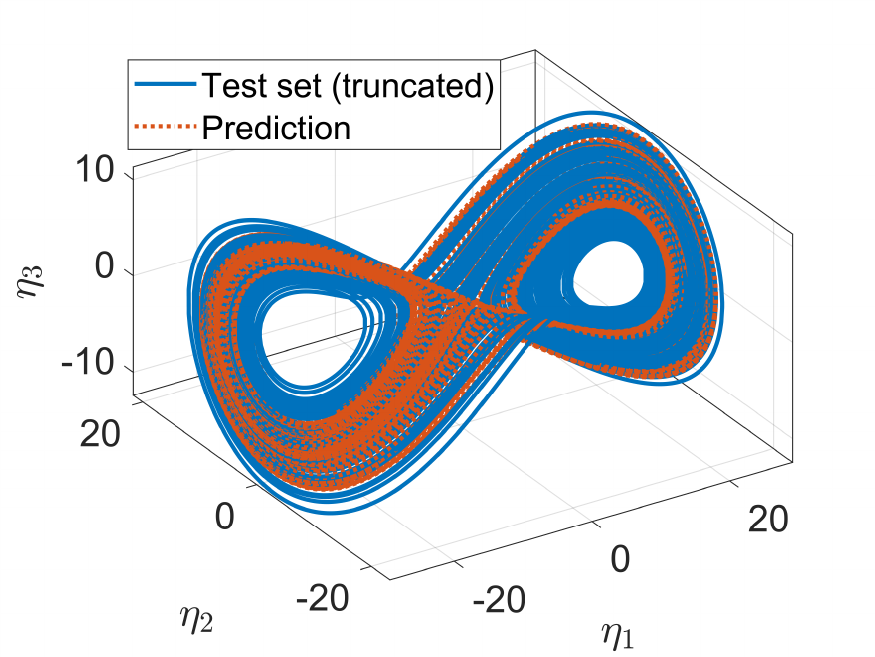}
    \caption{Reconstruction and prediction results of Lorenz attractor in 3D. The inertial manifold is reconstructed by observing only $x$ coordinate as a scalar time serie, delay embedded into 7D, and performing model reduction to 3D. Dynamics on SSM is represented by a fifth-order polynomial.}
    \label{fig:3D_ssm}
\end{figure}

In delayed coordinates, the origin is a fixed point corresponding to the origin of system (\ref{lorenz_eqa}). We therefore seek a 3D SSM attached to the origin in the 7-dimensional embedding space and determine the reduced dynamics on that SSM. By computing the normalized invariance error defined in Eq. (\ref{NMFE}) and identifying its minimum, we choose the order $\mathcal{K} = 3$ for the polynomial expansion for the SSM. We find that the 3D reduced dynamics on the SSM does admit a chaotic attractor, as shown in Fig. \ref{fig:3D_ssm}.

In the modal coordinates (\ref{modaldynamics}) that diagonalize the linear part of the SSM-reduced dynamics, we obtain the polynomial reduced-order model on the 3D SSM in the form 
\begin{equation} \label{Lorenz_RD}
    \begin{aligned}
        \dot{u} = & -11.0194u - 0.0001u^2 -0.0005 uv -0.0025 uw -0.0005 v^2  + 0.0026 w^2 + 0.0317 u^3 \\
        &  + 0.0622 u^2v -0.0011 u^2w + 0.0421 uv^2 +0.0102 uvw -0.0143 uw^2 +0.0058 v^3 \\
        &  + 0.0235 v^2w -0.0349 vw^2 + 0.0108 w^3 ,\\
        \dot{v} = & -2.3859v  +0.0036 uw +0.0002 v^2 + 0.0006 vw - 0.0025 w^2 -0.0273 u^3 - 0.0550 u^2v \\
        &  -0.0001 u^2w - 0.0391 uv^2 - 0.0096 uvw +0.0117 uw^2 -0.0059 v^3 - 0.0231 v^2w\\
        &  +0.0357 vw^2 - 0.0131 w^3 ,\\
        \dot{w} = & + 9.3871w - 0.0001 u^2 -0.0004 uv +0.0032 uw  + 0.0009 vw - 0.0016 w^2 -0.0144 u^3 \\
        &  - 0.0303 u^2v -0.0001 u^2w - 0.0233 uv^2 - 0.0054 uvw +0.0058 uw^2 -0.0041 v^3 \\
        &  - 0.0147 v^2w +0.0246 vw^2 - 0.0112 w^3 . \;
    \end{aligned}
\end{equation}

One might wonder what we have achieved here, given that system (\ref{Lorenz_RD}) is substantially more complicated than the original system (\ref{lorenz_eqa}) and has the same dimension as that system. One should not forget, however, that the input to our procedure was a single scalar observable and we used no information about the dimension or the form of the Lorenz system (\ref{lorenz_eqa}). The procedure we follow here is the same for any system that is only known from a scalar observable. This procedure will generally yield a smooth manifold, the SSM, that contains the attractor of the system. 


\begin{figure}[h!]
    \centering
    \begin{subfigure}{0.4\textwidth}
        \includegraphics[width=0.9\textwidth, height=2in]{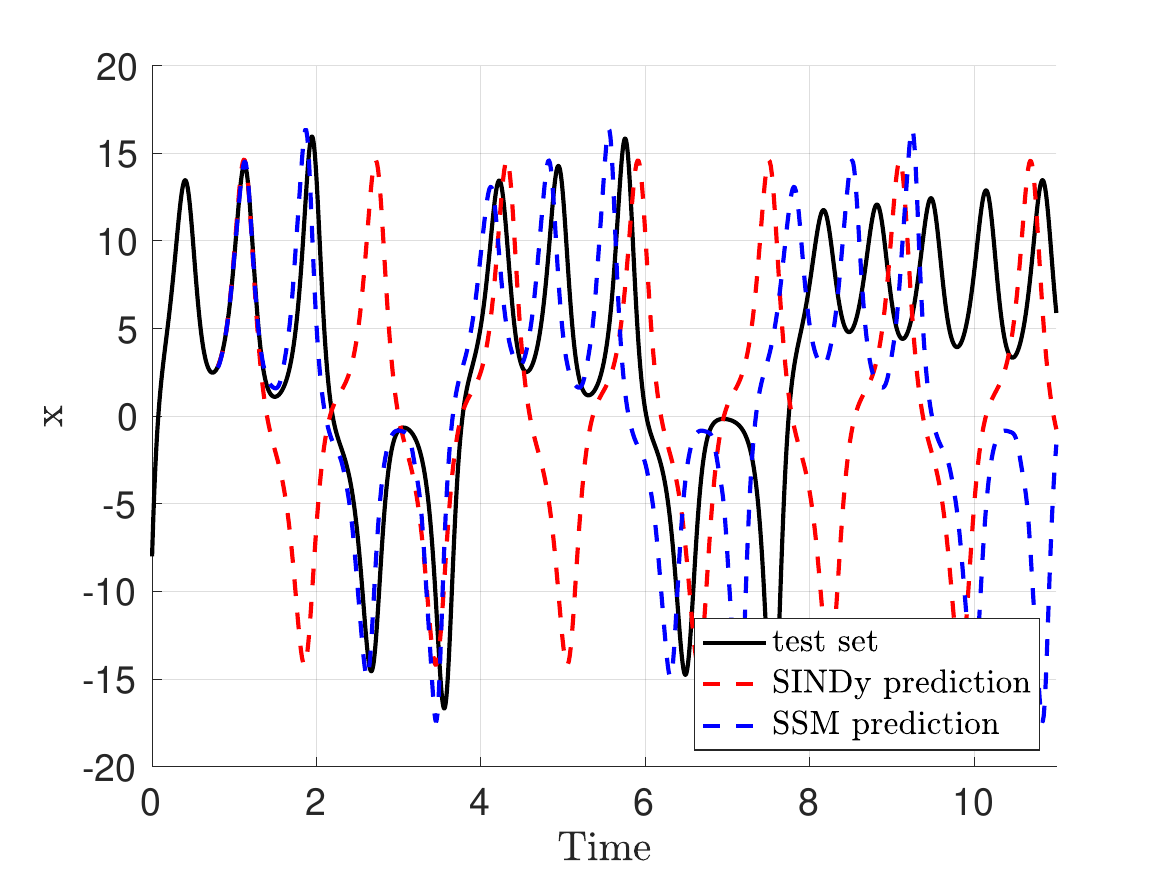}
        \caption{\label{fig:3D_pre}}
    \end{subfigure}
    \begin{subfigure}{0.4\textwidth}
        \includegraphics[width=0.9\textwidth, height=2in]{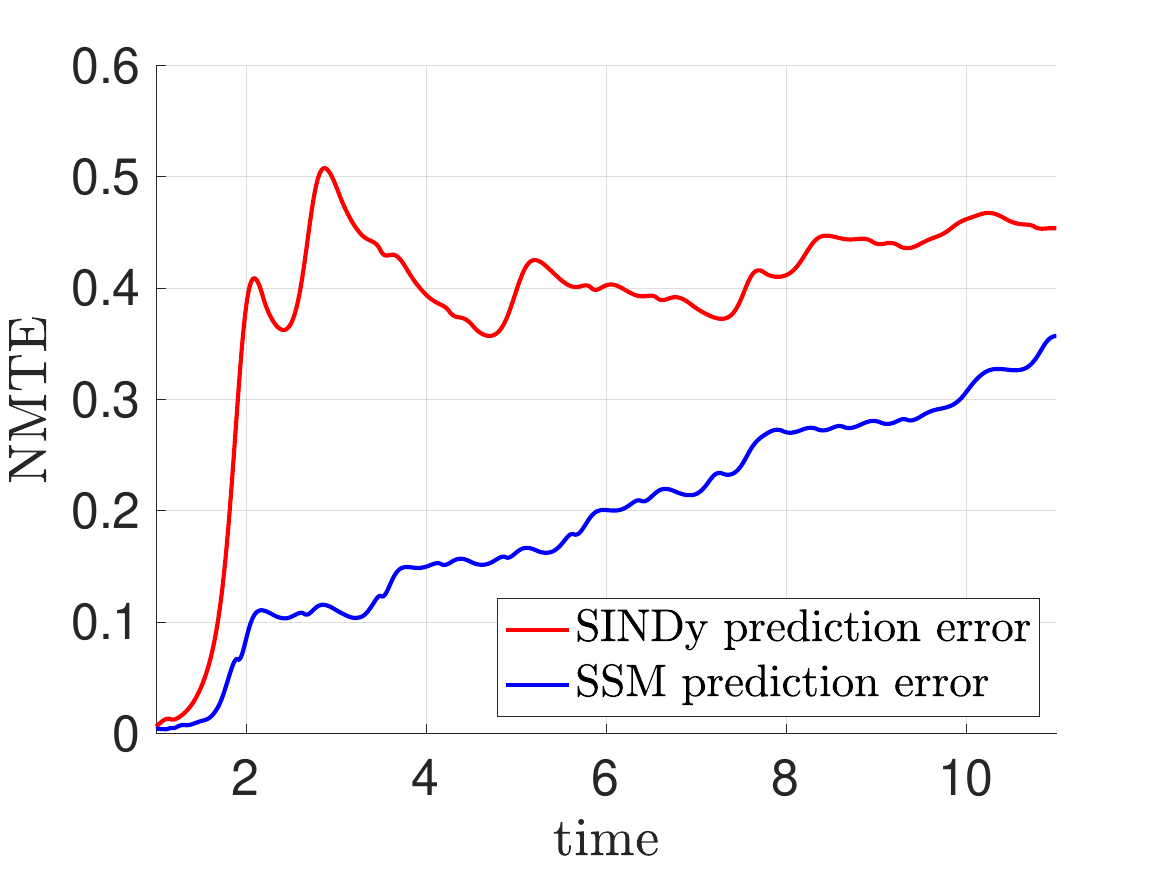}
        \caption{\label{fig:3D_NMTE}}
    \end{subfigure}
    \begin{subfigure}{0.4\textwidth}
        \includegraphics[width=0.9\textwidth, height=2in]{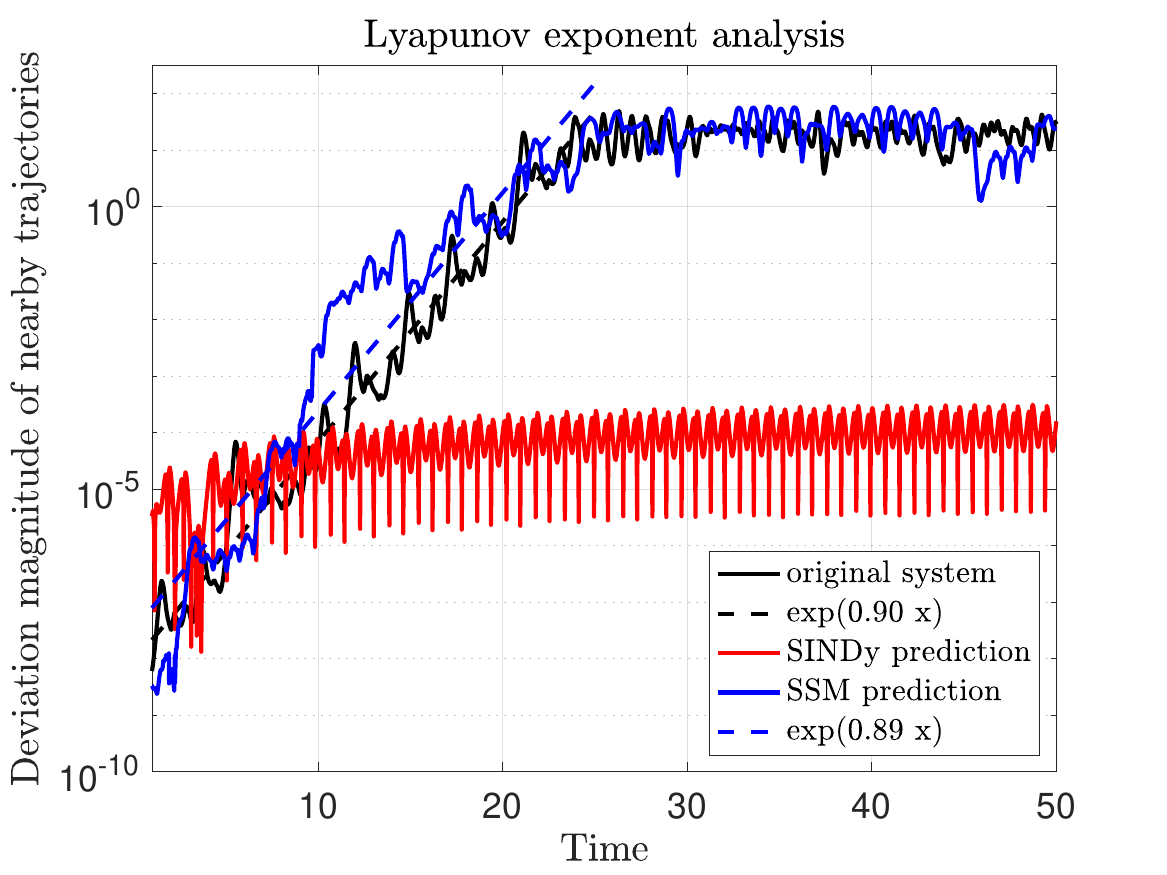}
        \caption{\label{fig:3D_le}}
    \end{subfigure}
    \caption{Time series forecasting results for delay-embedded Lorenz attractor, and its reconstructed Lyapunov exponent. Subplots (\subref{fig:3D_pre}) and (\subref{fig:3D_NMTE}) show more accurate predictions from SSM-reduction than from SINDy over approximately 2.97 Lyapunov times (NMTE less than 0.1). Subplot (\subref{fig:3D_le}) plots the separation of nearby trajectories against time. The red prediction made by the SINDy model has the maximal Lyapunov exponent 0.23 compared with the real MLE 0.90 of system (\ref{lorenz_eqa}). A blue linear fit to the SSM prediction confirms that trajectories separate exponentially, and the system is chaotic with the maximal Lyapunov exponent 0.89. At the same time, the red prediction made by SINDy suggests that the system is not chaotic. } 
    \label{fig:3D}
\end{figure}

With the SSM-reduced model, we can now make predictions for test trajectories not used in constructing this model, and compare it with predictions obtained from the SINDy algorithm. From the same training data, SINDy \cite{sindy} constructs a Hankel matrix by stacking 10 delayed time series as rows, and uses singular value decomposition (SVD) to select the first three dominant modes in these time-delay coordinates. Each column of the data is then trained with a different sparsification parameter ($\lambda _{\text{SINDy}} = 0.01, 0.2, 2$) to fit a third-order polynomial right-hand side of an ODE to the data. 

Figure \ref{fig:3D} shows the comparison of both short- and long-term predictions from SINDy and SSM reduction. Both methods are able to forecast for short times but our approach shows a more accurate result (Fig. \ref{fig:3D_pre} and Fig. \ref{fig:3D_NMTE}). We then perform a Lyapunov analysis by plotting the separation of two trajectories against time in a log plot in Fig. \ref{fig:3D_le}. The predictions given by \textsc{SSMLearn} admit a clear linear fit which indicates our model is chaotic, while the trajectories reconstructed by SINDy are not chaotic. Indeed, as noted in Section 4.5 in the Appendix of Ref. \cite{sindy}, although the skeleton of the Lorenz attractor is captured by SINDy, the attractor is reconstructed as a non-chaotic, quasi-periodic orbit. In contrast, the Lyapunov exponent from the SSM-based model (\ref{Lorenz_RD}) is $0.89$, which is very close to the real value of $0.90$ computed directly from the Lorenz system (\ref{lorenz_eqa}).

\subsection{A nine-dimensional Lorenz model}

We consider a 3D viscous fluid layer that is uniformly heated from below, and let $t$ and $\bm{r}$ be the time and space variables. The density field $\rho (\bm{r},t)$ describes the mass distribution of the fluid, and the physical state is determined by the pressure $\bm{p}(\bm{r},t)$, the temperature $\bm{T}(\bm{r},t)$ and the velocity field $\bm{v}(\bm{r},t)$. The dynamic behavior is governed by three nonlinear PDEs, namely the equation of continuity, the Navier-Stokes equation, and the equation of thermal conductivity:
\begin{equation}
    \begin{aligned}
    & \frac{\partial \bm{v}}{\partial t}+(\bm{v} \cdot \nabla )\bm{v} = \frac{\rho}{\rho _0} \bm{g} - \frac{1}{\rho _0}\nabla \bm{p} + \nu \nabla ^2 \bm{v} ,  \quad    \mathrm{div} \; \bm{v} = 0 , \\
    & \frac{\partial \bm{T}}{\partial t}+(\bm{v} \cdot \nabla ) \bm{T} = \chi \nabla ^2 \bm{T} ,     \;
    \end{aligned}
\end{equation}
where $\bm{g} = (0,0,-g)$ denotes the vector of gravity, $\rho _{0}$ the standard density, $\nu$ the kinematic viscosity and $\chi$ the thermal conductivity.

Lorenz applied a double Fourier expansion in deriving his 3D governing equations (\ref{lorenz_eqa}) in Ref. \cite{Lorenz}. Ref. \cite{9DLorenz} applied a similar approach but extended the double Fourier expansion to triple expansion. After truncation and selection of the leading-order-terms in the expansion, Ref. \cite{9DLorenz} obtains the nine-dimensional ODE 
\begin{equation} \label{9Dlorenz_eqa}
    \begin{aligned}
    \dot{C}_1 &= -\sigma b_1 C_1 -C_2 C_4 + b_4 C_4^2 + b_3 C_3 C_5 -\sigma b_2 C_7 , \\
    \dot{C}_2 &= -\sigma C_2 + C_1 C_4 -C_2 C_5 + C_4 C_5 -\frac{\sigma}{2}C_9    , \\
    \dot{C}_3 &= -\sigma b_1 C_3 +C_2 C_4 -b_4 C_2^2 -b_3 C_1 C_5 + \sigma b_2 C_8 ,\\
    \dot{C}_4 &= -\sigma C_4 + C_2 C_3 -C_2 C_5 + C_4 C_5 +\frac{\sigma}{2}C_9 ,\\
    \dot{C}_5 &= -\sigma b_5 C_5 + \frac{1}{2} C_2^2 -\frac{1}{2} C_4^2 ,\\
    \dot{C}_6 &= -b_6 C_6 + C_2 C_9 - C_4 C_9 , \\
    \dot{C}_7 &= -b_1 C_7 -r C_1 + 2 C_5 C_8 -C_4 C_9 , \\
    \dot{C}_8 &= -b_1 C_8 +r C_3 - 2 C_5 C_7 +C_2 C_9 , \\
    \dot{C}_9 &= -C_9 -r C_2 + r C_4 -2 C_2 C_6 + 2 C_4 C_6 + C_4 C_7 - C_2 C_8 , \;
    \end{aligned}
\end{equation}
where $\sigma$ is the Prandtl number, $r$ is the Rayleigh number, and the constant parameters $b_i$ represent a measure of the geometry of the square cell.

Depending on $r$ and $\sigma$, the dynamics can be stationary, periodic, chaotic or hyperchaotic, with the latter referring to the case wherein multiple positive Lyapunov exponents exist. For certain parameters, four symmetric low-dimensional chaotic attractors are observed in system (\ref{9Dlorenz_eqa}). From \cite{9DLorenz}, we choose the system parameters as
\begin{equation}
    \begin{aligned}
        \sigma & = \frac{1}{2}, \; r  =14.2, \;a = \frac{1}{2} , \\
        b_1 & = \frac{4(1+a^2)}{1+2a^2}, \; b_2 = \frac{1+2a^2}{2(1+a^2)}, \; b_3 = \frac{2(1-a^2)}{1+a^2} , \\
        b_4 &= \frac{a^2}{1+a^2}, \; b_5 = \frac{8a^2}{1+2a^2}, \; b_6 = \frac{4}{1+2a^2} .\;
    \end{aligned}
\end{equation}

\begin{table}[h!]
\centering
\begin{tabular}{ccccccccc}
\toprule
$\lambda_1$ & $\lambda_2$ & $\lambda_3$ & $\lambda_4$ & $\lambda_5$ & $\lambda_6$ & $\lambda_7$ & $\lambda_8$ & $\lambda_9 $\\ \hline
1.9263                     & -0.2741                    & -0.2741                    & -0.5000                    & -0.6667                    & -2.6667                    & -3.4263                    & -4.7259                    & -4.7259    \\ \bottomrule               
\end{tabular}
\caption{Eigenvalues of the 9D Lorenz attarctor.}
\label{9D_eigenvalues}
\end{table}

The origin is a fixed point of system (\ref{9Dlorenz_eqa}) with the eigenvalues listed in Table (\ref{9D_eigenvalues}). As in our previous example, no information about these eigenvalues will be used in our data-driven SSM reduction; we only list these eigenvalues for reference. Note that the 1:1 resonances between $\lambda_2$ and $\lambda_3$ and between $\lambda_8$ and $\lambda_9$ do not violate the non-resonance condition (\ref{nonResonance}). There is also a 5th-order resonance between $\lambda_4$, $\lambda_5$ and $\lambda_6$, which violate Eq. (\ref{nonResonance}). However, as earlier work on SSMs shows \cite{SSM,SSMLearn}, a polynomial expansion for an SSM, $\mathcal{W}(E)$, can still be computed as long as the linearized spectrum within $E$ does not contain $\lambda_4$, $\lambda_5$ and $\lambda_6$. We will choose $E$ in this fashion.

\begin{figure}[h]
    \centering
    \includegraphics[width=0.5\textwidth]{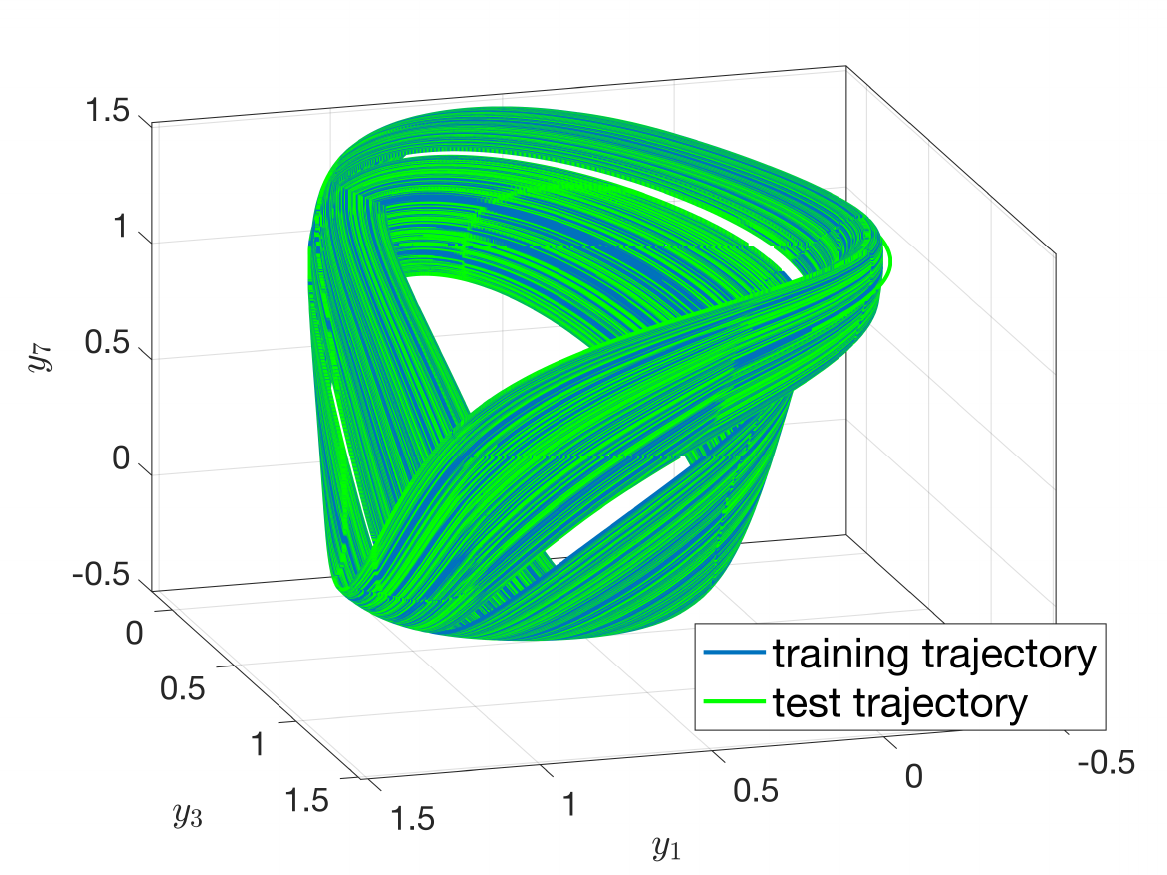}
    \caption{The chaotic attractor of the 9D Lorenz model (\ref{9Dlorenz_eqa}) in the 7D delay-embedded space, projected onto the coordinates $y_1$, $y_3$ and $y_7$.}
    \label{fig:9D_ss}
\end{figure}

As in the 3D Lorenz example, we observe only the $C_1$ coordinate of system (\ref{9Dlorenz_eqa}) as a scalar time series and delay-embed the signal. We generate two random initial conditions (one for training and the other for testing) within the unit sphere of the origin in the 9D phase space and observe the first coordinate $C_1$. We then select and truncate the trajectories to make them lie close to one of the four symmetric strange attractors. Assuming no knowledge of the original system or the dimension of the chaotic attractor, we again use FNN to estimate the dimension $d$ of the SSM, and pick $d = 3$ according to the false neighbor percentage in Table \ref{9D_dim}. By the Takens embedding theorem, we delay-embed this 3D SSM into a $\rho = 2 \times d +1  = 7$ dimensional space. The training and testing trajectories in this seven-dimensional delay-embedded space are shown in Fig. \ref{fig:9D_ss}.

\begin{table}[h]
\centering
\begin{tabular}{cccccccc}
\toprule
SSM dimension ($d$)           & 2 & 3 & 3 & 3 & 3 & 3 & 3 \\ \midrule
FNN percentage           & 41.6\%  & 0.0\%  & 0.0\%  & 0.0\% & 0.0\% & 0.0\% & 0.0\% \\ \midrule
Delay embedding dimension  ($\rho$)    & 5 & 7 & 7 & 7 & 7 & 7 & 7 \\ \midrule
Manifold expansion order ($\mathcal{K}$) & - & 2 & 3 & 4 & 5 & 6 & 7 \\ \midrule
Invariance Error         & - & 2.24\%  & 0.59\%  & 0.22\% & 0.09\% & 0.05\% & 0.04 \% \\ \bottomrule
\end{tabular}
\caption{Different choices of SSM dimension, false nearest neighbor percentage (FNN), order of manifold parametrization, and the corresponding invariance error of the delay-embedded 9D Lorenz model. Due to the high FNN percentage for $d=2$, we did not compute an SSM approximation for that case.}
\label{9D_dim}
\end{table}

Using \texttt{SSMLearn}, we approximate a 3D SSM attached to the origin of the 7D delay-embedded space. We minimize the invariance error over different choices of SSM dimension and polynomial approximation orders. As shown in Table \ref{9D_dim}, these considerations lead us to construct a 3D SSM via a sixth-order polynomial expansion. This construct produces an invariance error of order $0.05\%$, which is low enough to guarantee high accuracy for our SSM-reduced model.

\begin{figure}[h!]
    \centering
    \begin{subfigure}{0.6\textwidth}
        \includegraphics[width=\textwidth, height=2in]{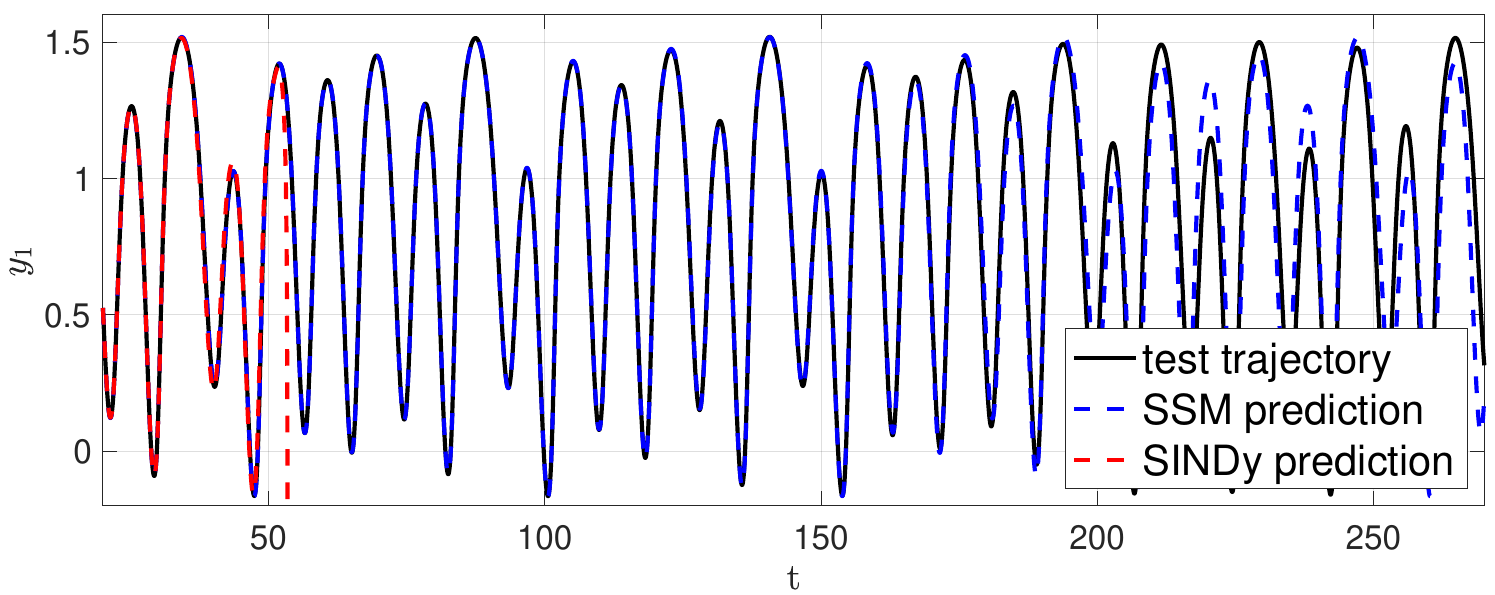}
        \caption{\label{fig:9D_pre}}
    \end{subfigure}
    \begin{subfigure}{0.39\textwidth}
        \includegraphics[width=\textwidth, height=2in]{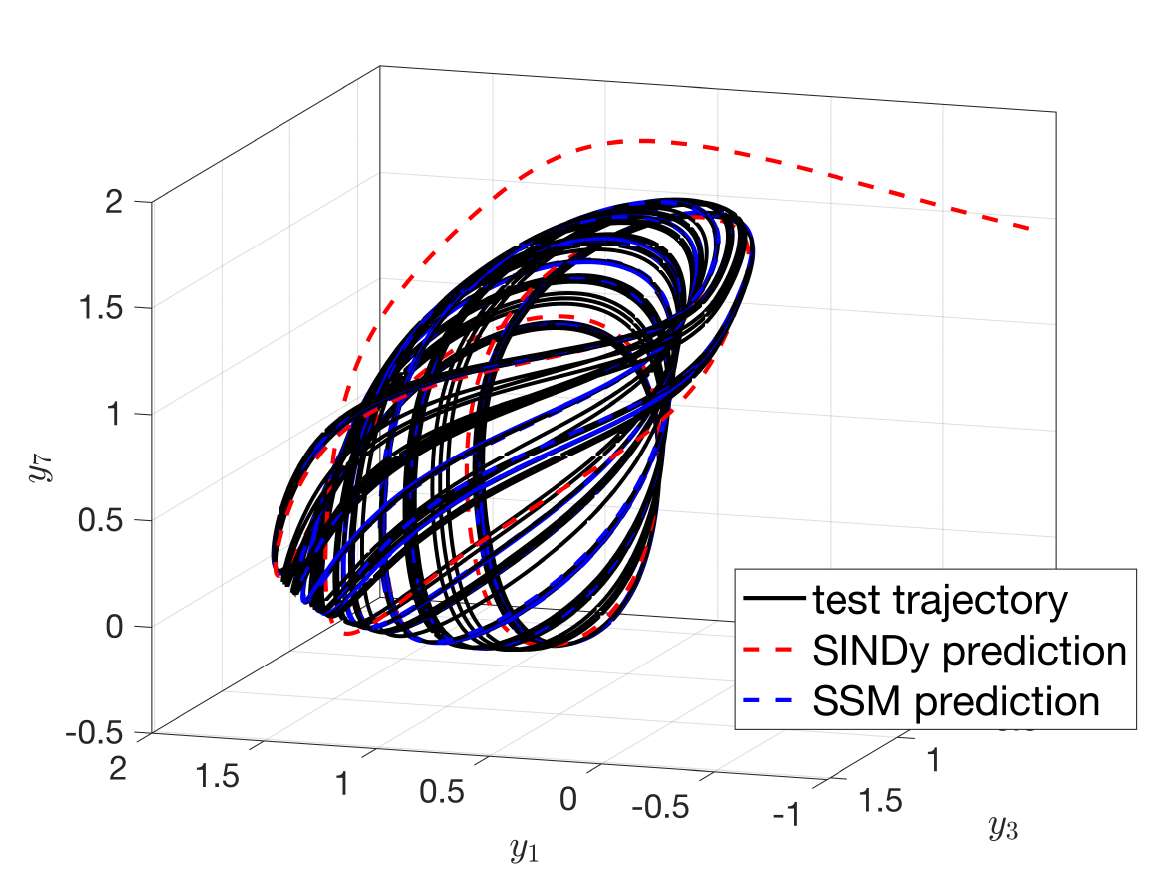}
        \caption{\label{fig:9D_pre2}}
    \end{subfigure}
    \caption{Trajectory prediction in the time domain (plot (\subref{fig:9D_pre})) and in the delay-embedding space (plot (\subref{fig:9D_pre2})) from the reduced models for an initial condition not contained in the training data. The length of the SSM prediction (blue) interval is approximately 5.94 Lyapunov times, while the SINDy model (red) fails for trajectory integration at around 1.82 Lyapunov times. We obtain the prediction of the SSM-reduced model by making predictions on 3D reduced coordinates using the kNN method, then projecting the trajectory back in 7D delay-embedded space. The training set contains $1.5$ million points.} 
    \label{fig:9D_prediction}
\end{figure}

In the reduced coordinates, we find that a polynomial representation of the SSM-reduced dynamics fails to capture the complex dynamics with sufficient accuracy. Instead, we turn to the $k$-th nearest neighbor method discussed in Section \ref{sec:prediction}. The number of neighbors in kNN algorithm is often chosen to be larger than the data dimension \cite{knnFarmer}. Here, we select $k=4$ and perform the kNN method on a training set containing $1.5$ million data points. Note here that the training data for SSM fitting ($10^5$ points) is in general much less than the number of points required for reduced dynamics fitting. This is because the kNN prediction algorithm requires data to fully cover the attractor in its phase space. Fig. \ref{fig:9D_pre} shows the prediction results from this nearest neighbor method. Since chaotic dynamics have sensitive dependence on initial conditions, the forecast is only accurate for 5.94 Lyapunov times.

We attempted to perform the SINDy algorithm on the same training data to identify a seven-dimensional SINDy model and compare it to SSM-reduced model prediction results. We tested several combinations of different model orders (from two to six) and sparsity parameters ($\lambda _{\text{SINDy}} = 0.0001$ to $0.1$).  We also tried to follow the SINDy tutorial for \cite{Kaptanoglu2022} to tune the hyperparameter $\lambda _{\text{SINDy}}$ by minimizing RMSE of predicted trajectories for different choices of $\lambda _{\text{SINDy}}$. However, the resulting SINDy models do not contain a global attractor in general, and their solutions tend to blow up after some time of integration. As an example, Fig. \ref{fig:9D_prediction} shows results from a third-order SINDy model with $\lambda _{\text{SINDy}} = 0.001$. One could also try to set different sparsity parameters for each coordinate, as for the delay-embedded Lorenz example in \cite{sindy}, but we did not pursue that exploration here.

\begin{figure}[h!]
    \centering
    \begin{subfigure}{0.49\textwidth}
        \includegraphics[width=\textwidth, height=2in]{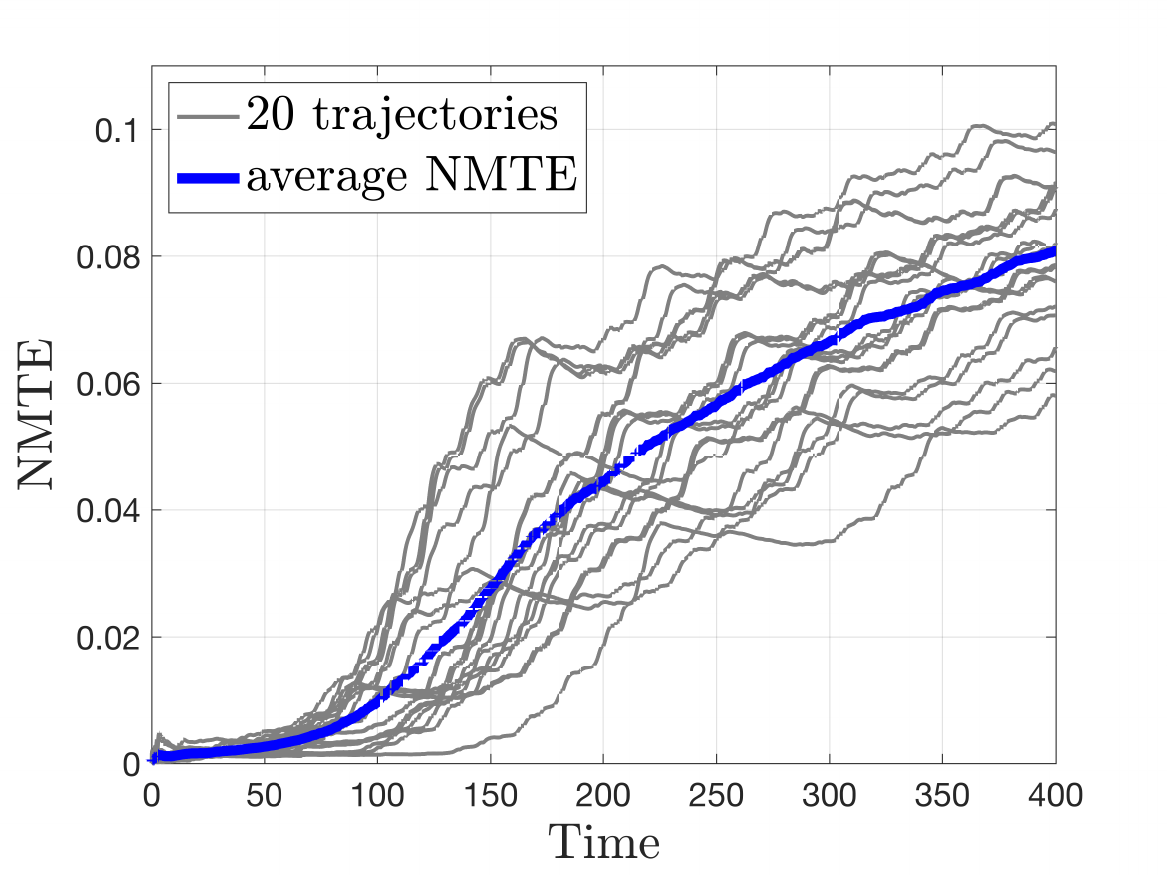}
        \caption{\label{fig:9D_NMTE}}
    \end{subfigure}
    \begin{subfigure}{0.49\textwidth}
        \includegraphics[width=\textwidth, height=2in]{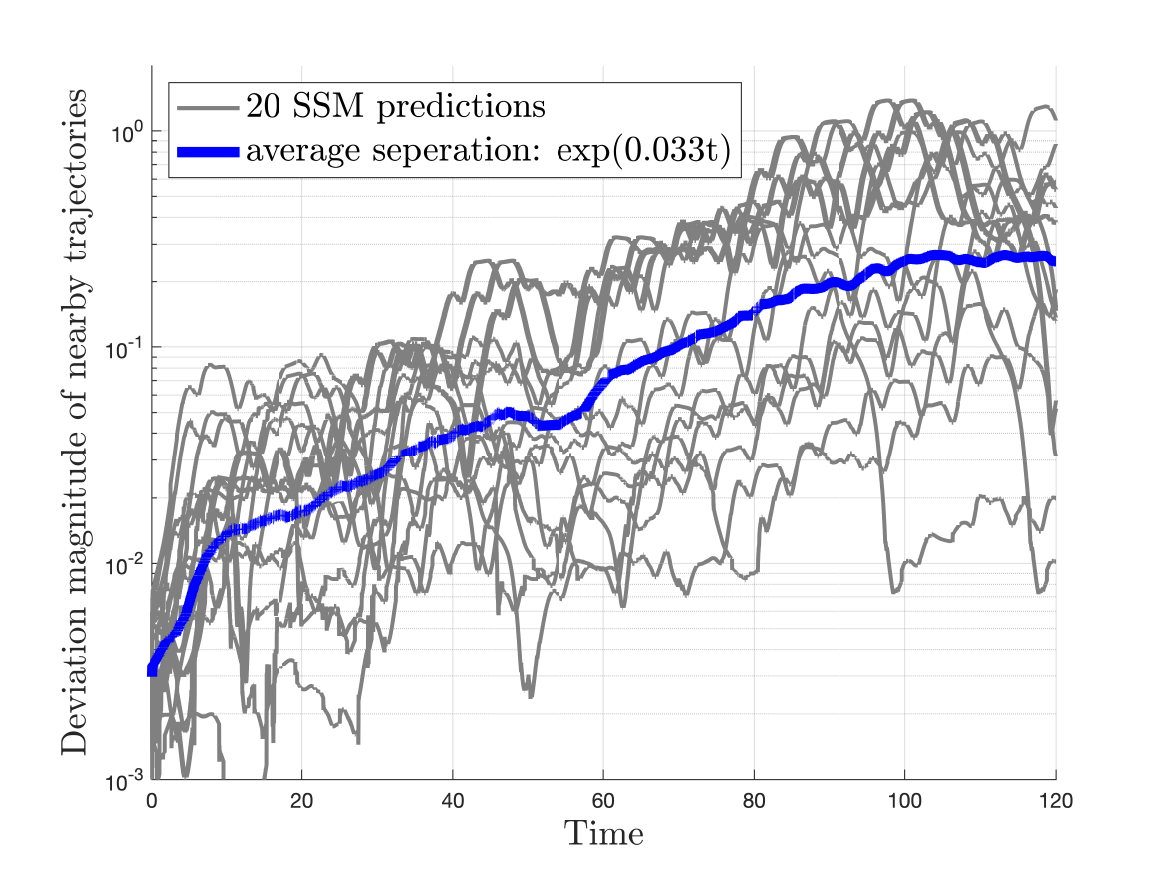}
        \caption{\label{fig:9D_le}}
    \end{subfigure}
    \caption{Trajectory prediction error (\subref{fig:9D_NMTE}) and the Lyapunov exponent analysis (\subref{fig:9D_le}) on test data of the 3D SSM-reduced model. We compute the time average of the NMTE and the nearby trajectory separation rate over 200 test trajectories, but only plot 20 of them in this figure. The average trajectory separation rate of system (\ref{9Dlorenz_eqa}) is 0.032, which is closely approximated by the mean Lyapunov exponent 0.033 computed from subplot (\subref{fig:9D_le}) of the SSM-reduced model.} 
    \label{fig:9D_20}
\end{figure}


To test the statistical properties of the SSM-reduced model, we further compute the NMTE and Lyapunov exponent of 200 pieces of the test trajectory whose initial conditions are not contained in the training data. After adding a small perturbation, we plot the deviation of the original and the predicted model trajectory against time on a log scale. For 200 trajectories of the original system (\ref{9Dlorenz_eqa}), this analysis yields the Lyapunov exponent 0.032 as the mean of the observed separation exponents. Remarkably, this exponent is closely approximated by the value 0.033 yielded by our SSM-reduced model in Fig. \ref{fig:9D_le}.



\subsection{Duffing oscillator chain} \label{sec:duffing}

\begin{figure}[h!]
    \centering
    \includegraphics[width=0.7\textwidth]{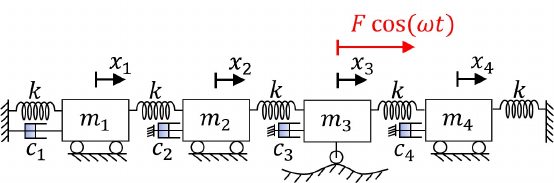}
    \caption{Forced harmonic oscillator chain with a single Duffing oscillator added. We apply periodic forcing to this Duffing oscillator. }
    \label{fig:FD}
\end{figure}

We now consider the classic forced-damped Duffing oscillator (\cite{HamelGeorgDI}) given by
\begin{equation}
    \ddot{x} + \delta \dot{x} + \alpha x + \beta x^3 = \gamma \cos{\omega t}, 
\end{equation}
where the $\delta$, $\alpha$, $\beta$ and $\gamma$ are the coefficients for damping, linear stiffness, nonlinear stiffness, and forcing amplitude, respectively. We use this nonlinear oscillator as the third element in an otherwise linear chain of oscillators, as shown in Fig. \ref{fig:FD}. We denote the displacement of the $i$-th degree of freedom by the coordinate $x_i$. We then apply periodic forcing along $x_3$ as illustrated in Fig. \ref{fig:FD}.

We choose the third mass as $m_3 = 1$ kg and the remaining mass as $m_1 = m_2 = m_4 = 0.1$ kg each; the third linear damping coefficient attached to the ground is $c_3 = 0.1$ N$\cdot$s/m; the other damping coefficients are $c_1 = c_2 = c_3 = c_4 = 0.75$ N$\cdot$s/m. The linear spring stiffness is $k = 1$ N/m and the third nonlinear spring has the force relation $F_3(x_3) = k_3 x_3 + \beta x_3^3$, where $k_3 = 3$ N/m and $\beta = -0.25$ N/$\text{m}^3$. The periodic forcing that acts on the Duffing oscillator has an amplitude of 2.25 N and a frequency of 2 rad/s. The resulting equations of motion is given by a system of second-order ODEs of the form
\begin{equation} \label{sys_duffing}
    \begin{aligned}
    0.1\ddot{x}_{1}+0.75\dot{x}_{1} +2x_{1}-x_{2} &= 0 ,\\
    0.1\ddot{x}_{2}+0.75\dot{x}_{2} -x_{1}+2x_{2}-x_{3} &= 0 ,\\
    \ddot{x}_{3}+0.1\dot{x}_{3} -x_{2}-x_{3} - x_{4} + 0.25x_{3}^{3}&= 2.25 \cos{2t}, \\
    0.1\ddot{x}_{4}+0.75\dot{x}_{4} -x_{3}+3x_{4} &= 0. \;
    \end{aligned}
\end{equation}

When written in first-order form, these equations define an eight-dimensional, non-autonomous dynamical system. For the autonomous part of the dynamics, the origin is an unstable fixed point. This unstable equilibrium has eight eigenvalues, as listed in Table \ref{duffing_eigenvalues}. Again, these are listed only for reference and will not be used in our SSM construction.

\begin{table}[h]
\centering
\begin{tabular}{ccccc}
\toprule
$\lambda_{1}$ & $\lambda_2$ & $\lambda_{3,4}$ & $\lambda_{5,6}$ & $\lambda_{7,8}$ \\ \hline
1.2204  & -5.8047 & $-1.5713\pm0.6269i$ & $-3.6865\pm4.0004i$  & $-3.7500\pm3.9922i$   \\       \bottomrule               
\end{tabular}
\caption{Eigenvalues of the autonomous part of the Duffing oscillator chain.}
\label{duffing_eigenvalues}
\end{table}

\begin{figure}[h]
    \centering
    \includegraphics[width=0.5\textwidth]{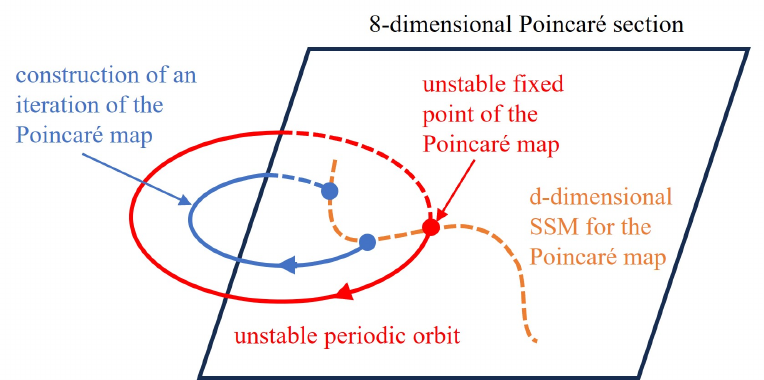}
    \caption{SSM construction for the Poincaré map of system (\ref{sys_duffing}) in the phase space.}
    \label{fig:SSM_2D}
\end{figure}

When we add the periodic forcing term to the autonomous part of the equation, the unstable fixed point at the origin bifurcates into an unstable periodic orbit. As long as the amplitude of the time-dependent term in system (\ref{sys_duffing}) is moderate, the corresponding $T = \pi$-periodic Poincaré map will also have an unstable fixed point. This hyperbolic fixed point is the anchor point for our SSM computation. Our strategy will be to locate an SSM for this map that contains an observed chaotic attractor of the system. As in our previous examples, this SSM will act as an inertial manifold whose restricted dynamics provide a reduced model for the dynamics on and near the attractor. The corresponding geometry is illustrated in Fig. \ref{fig:SSM_2D}. 

\begin{figure}[h]
    \centering
    \begin{subfigure}{0.4\textwidth}
        \includegraphics[width=0.9\textwidth, height=2in]{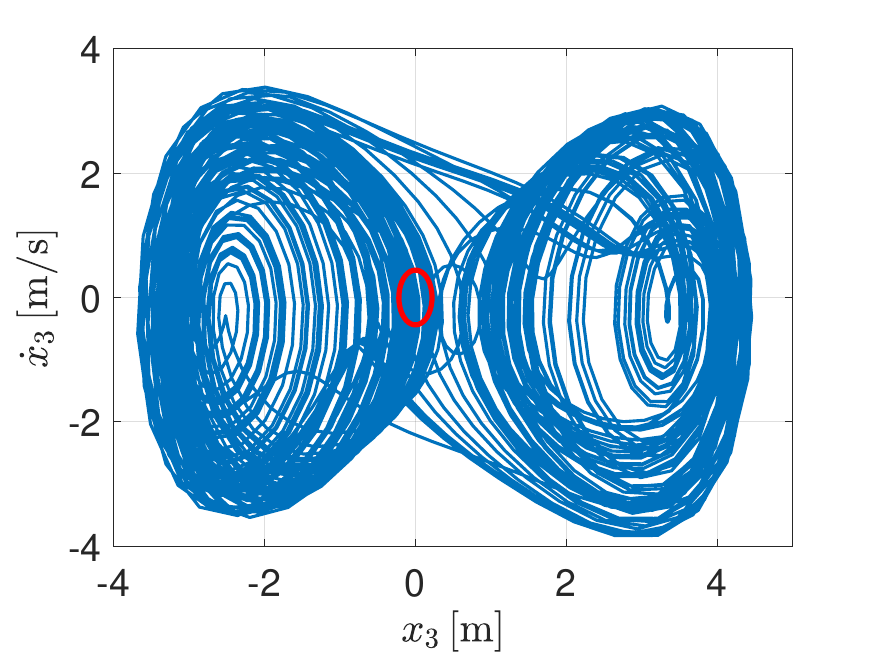}
        \caption{\label{fig:FD_flow}}
    \end{subfigure}
    \begin{subfigure}{0.4\textwidth}
        \includegraphics[width=0.9\textwidth, height=2in]{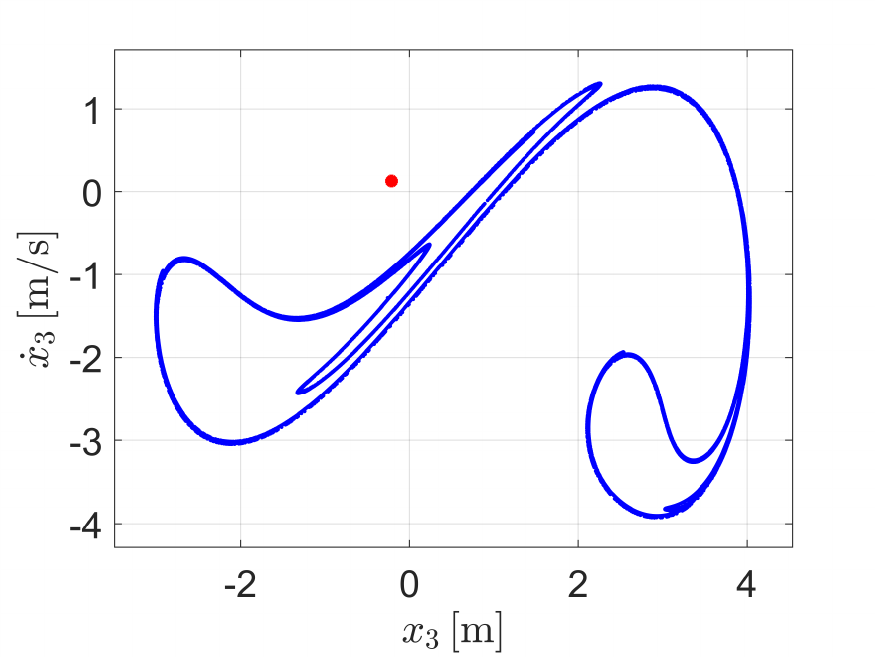}
        \caption{\label{fig:FD_map}}
    \end{subfigure}
    \caption{(\subref{fig:FD_flow}) Chaotic trajectories (blue) and a saddle-type periodic orbit of system (\ref{sys_duffing}) projected onto the $(x_3, \dot{x}_3)$ plane.
    (\subref{fig:FD_map}) The same trajectories and periodic orbits, as they appear as a fixed point and a trajectory on the attractor of the Poincaré map, projected again to the $(x_3, \dot{x}_3)$ plane.} 
    \label{fig:FD_flowmap}
\end{figure}

For the eight-dimensional chaotic system (\ref{sys_duffing}), we consider all eight degrees of freedom as observables. We select seven initial conditions (six for training and one for testing) randomly within the unit sphere around the origin, and generate the corresponding trajectories for $10^5$ times the forcing period. We track the evolution of the origin under increasing forcing amplitude into a hyperbolic periodic orbit using the numerical continuation package \texttt{COCO} \cite{coco1}. In Fig. \ref{fig:FD_flow}, we plot the coordinates $x_{3}$ and $\dot{x}_{3}$ for the original continuous dynamics in blue, and the unstable periodic orbit computed by \texttt{COCO} in red. When we pass to the Poincaré map, the continuous orbits become discrete dots in the Poincaré section as in Fig. \ref{fig:FD_map}. Based on the theory of mixed-mode SSMs \cite{mixedmodeSSM}, \texttt{SSMLearn} computes a polynomial expansion for an SSM attached to the hyperbolic fixed point (red dot in Fig. \ref{fig:FD_map}). Note that our construct relies on the existence of the hyperbolic fixed point, which we compute using \texttt{COCO} from system (\ref{sys_duffing}). However, even when the exact location of this anchor point is not known, \texttt{SSMLearn} can still provide an approximation for the SSM automatically by regressing a manifold (with the addition of constant terms in its parameterization) to the available data. The location of the fixed point itself will then follow as a prediction from the SSM-reduced dynamics.

As in the previous examples, we use the invariance error as a criterion for determining the correct SSM dimension. Table \ref{FD_dim} shows different invariance errors obtained for different choices of SSM dimension and expansion order. Based on these results, we choose a 3D SSM construction via a 7th-order polynomial expansion.

\begin{table}[h]
\centering
\begin{tabular}{c|ccc|ccc}
\hline
SSM dimension  ($d$)   & \multicolumn{3}{c|}{2D SSM} & \multicolumn{3}{c}{3D SSM} \\ \hline
SSM expansion orders $(\mathcal{K})$ & 3     & 5     & 7     & 3     & 5     & 7    \\ \hline
Invariance Error         & 8.26\%  & 4.24\%  & 1.89\%  & 0.061\%  & 0.005\%  & 0.0001\% \\ \hline
\end{tabular}
 \caption{Different choices of SSM dimension, manifold expansion order and the corresponding invariance error for the oscillator chain (\ref{sys_duffing}).}
 \label{FD_dim}
\end{table}

\begin{figure}[h]
    \centering
    \includegraphics[width=1\textwidth]{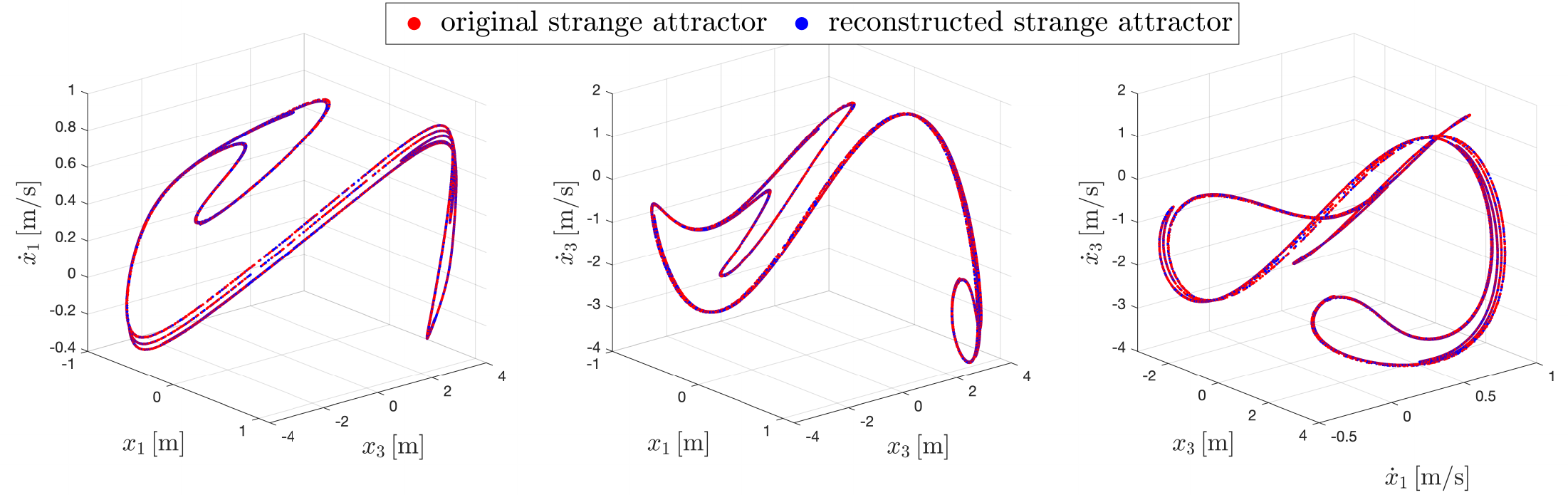}
    \caption{Predictions for the strange attractor of the forced Duffing oscillator chain in the 8-dimensional full state space, using the kNN prediction method on a three-dimensional SSM. Here $k=6$ and the training set contains $6 \times 10^{5}$ points.}
    \label{fig:FD_pre}
\end{figure} 

\begin{figure}[h]
    \centering
    \includegraphics[width=1\textwidth]{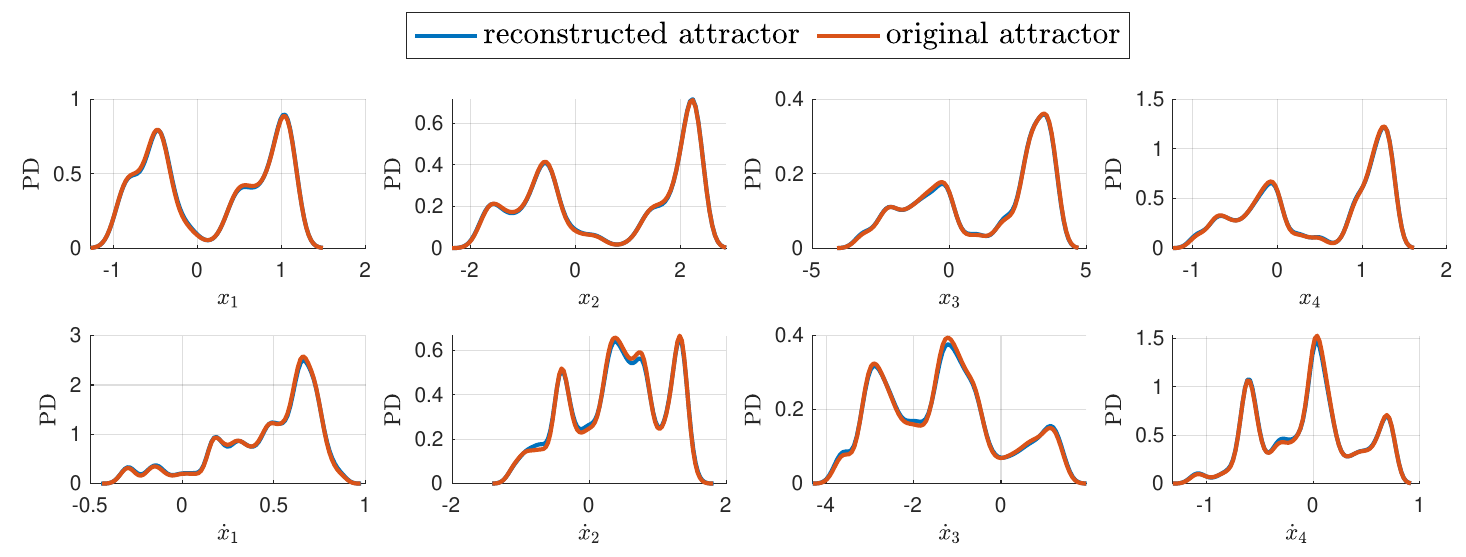}
    \caption{Probability distribution for the chaotic attractor of the full system (\ref{sys_duffing}) (red) and for the attractor of the 3D, SSM-reduced model in various sets of coordinates (blue). }
    \label{fig:FD_pd}
\end{figure}

After obtaining reduced coordinates on the 3D SSM, we again use the kNN algorithm discussed in Section \ref{sec:prediction} to reconstruct the chaotic attractor of the system on this SSM. We set $k=6$ and use a training set with the total size of $6 \times 10^{5}$ points. The prediction is generated by the 3D reduced model, then lifted to the eight-dimensional state space via the SSM parametrization. In Fig. \ref{fig:FD_pre}, we visualize the test data (red) and its prediction (blue) by showing them in various combinations of three of the eight coordinates. From Figure \ref{fig:FD_pre}, we see that the red dots and the blue dots coincide in the 8D state space, and the shape of the chaotic attractor is accurately reconstructed. To verify the accurate reproduction of long-term statistical features of the SSM-reduced model, we further plot the probability distribution of the original and the reconstructed dynamics in all of the coordinates separately. From Fig. \ref{fig:FD_pd}, the statistical properties of the reduced-order model are very close to those of the original attractor.

We also attempted to apply the map version of SINDy \cite{sindy_map} to the same eight-dimensional training data (six trajectories) for maps. However, after testing different combinations of various values of polynomial order (from three to six) and sparsity parameter ($\lambda _{\text{SINDy}} = 0.001 $ to $0.05$), we failed to identify a feasible model using SINDy and the prediction results invariably blew up after some time. This is likely due to the high complexity of the Poincaré map dynamics. 

\subsection{Delay-embedded Rössler attractor} \label{sec:rossler}

In our previous examples, we considered mixed-mode SSMs tangent to spectral subspaces with only real eigenvalues (i.e., $q = 0$ in equation (\ref{eqa_q})). Here, we turn to an example of an SSM tangent to an eigenspace with complex conjugate eigenvalues, the Rössler attractor system \cite{Rssler1976}
\begin{equation} \label{ross_eqa}
    \begin{aligned}
    \dot{x} &= -y-z ,\\
    \dot{y} &= x+ay  ,\\
    \dot{z} &= b+z(x-c).  \;
    \end{aligned}
\end{equation}

This system has two unstable fixed points at 
$$\left(\begin{array}{c}x_0\\ y_0\\ z_0\end{array}\right) = \left(\begin{array}{c}\frac{1}{2}(c \pm \sqrt{c^2-4ab}) \\ \frac{1}{2a}(-c \mp \sqrt{c^2-4ab}) \\ \frac{1}{2a}(c \pm \sqrt{c^2-4ab})\end{array} \right).$$
One of these fixed points lies in the center of the Rössler attractor loop and the other one lies relatively far from the attractor. We choose the classical parameter values $a = b = 0.2, c = 5.7$ and focus on the unstable equilibrium $(x_0, y_0, z_0) = (0.007, -0.035, 0.035)$ located in the center. We build the mixed-mode SSM attached to this fixed point that has one real eigenvalue $-5.687$ and a complex conjugate pair $0.097 \pm 0.995i$. 

\begin{figure}[h!]
    \centering
    \begin{subfigure}{0.4\textwidth}
        \includegraphics[width=0.9\textwidth, height=2in]{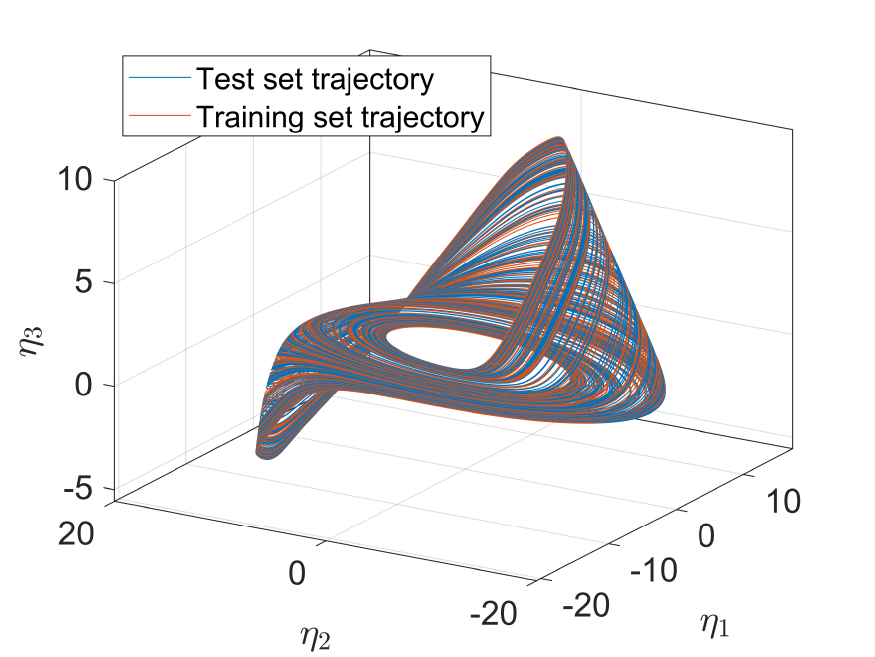}
        \caption{\label{fig:ross_rd}}
    \end{subfigure}
    \begin{subfigure}{0.4\textwidth}
        \includegraphics[width=0.9\textwidth, height=2in]{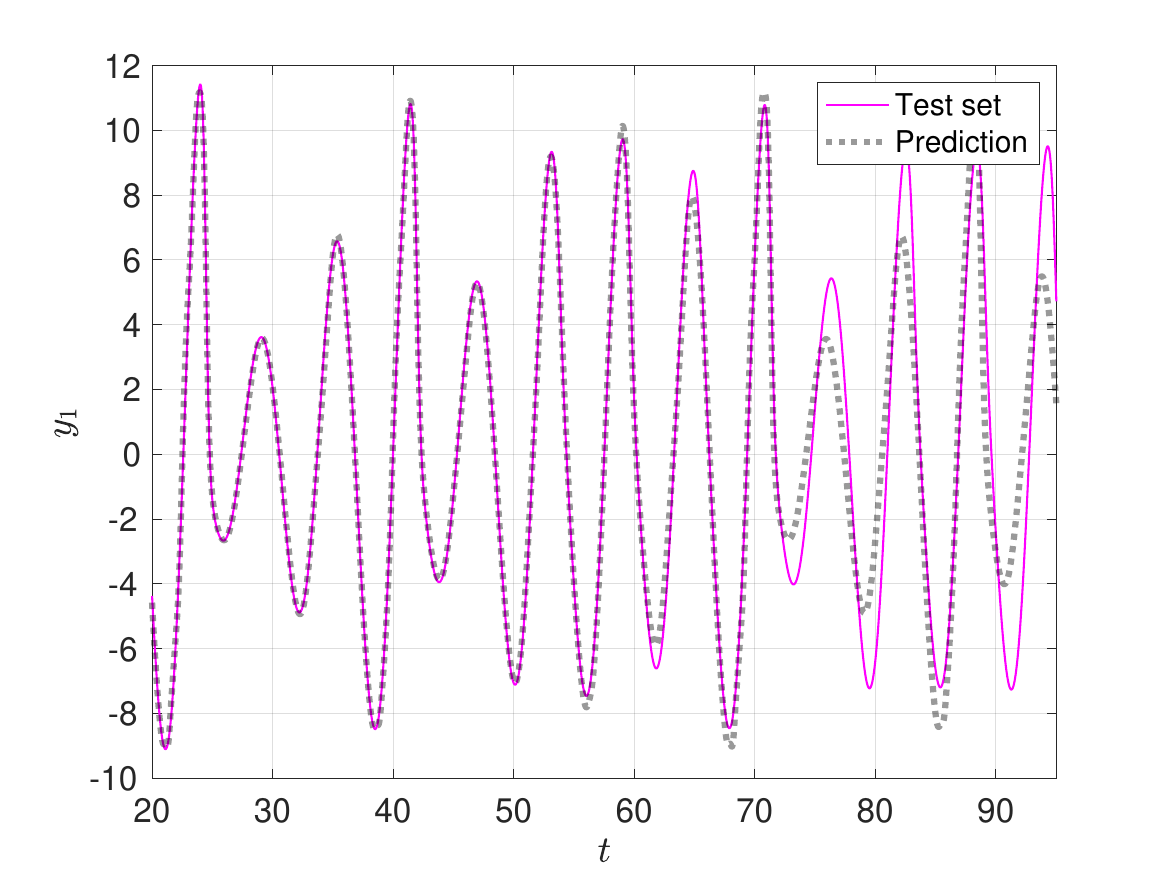}
        \caption{\label{fig:ross_pre}}
    \end{subfigure}
    \begin{subfigure}{0.4\textwidth}
        \includegraphics[width=0.9\textwidth, height=2in]{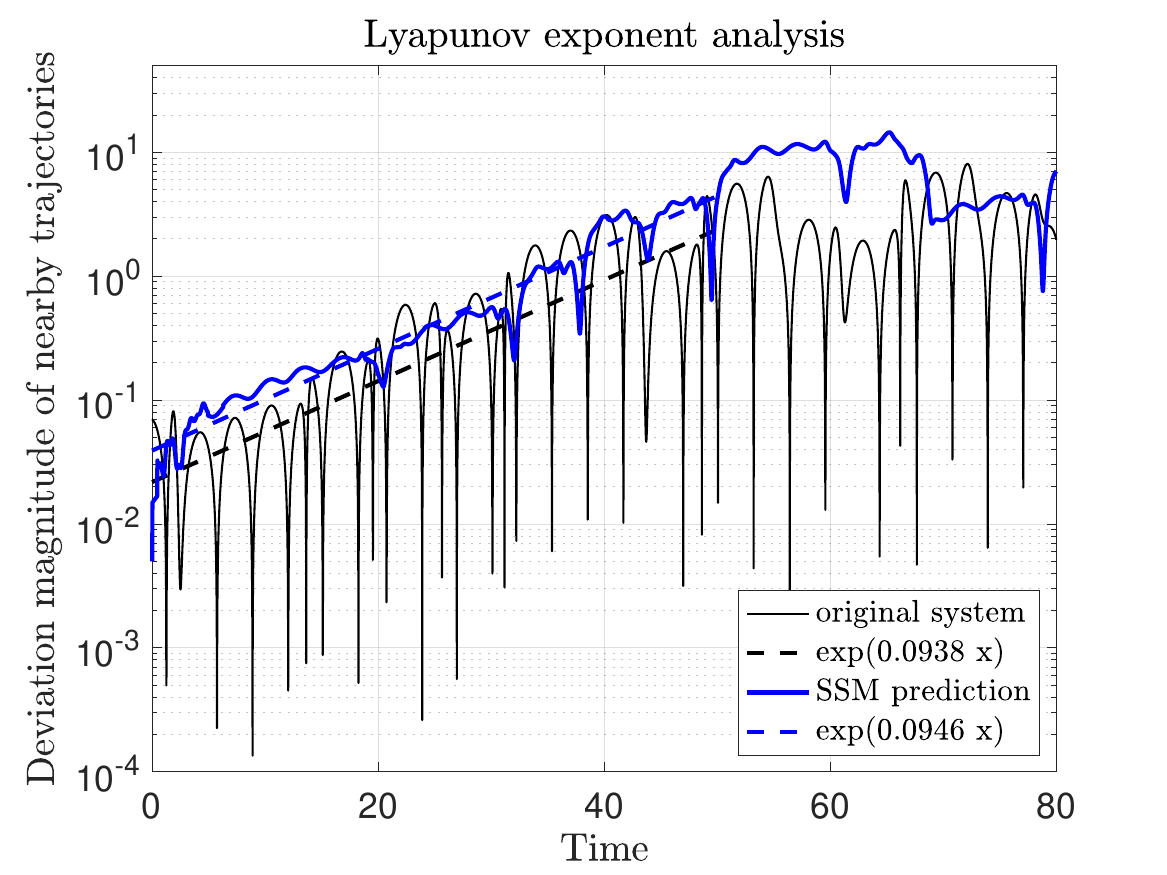}
        \caption{\label{fig:ross_le}}
    \end{subfigure}
    \caption{The delay-embedded Rössler attractor in reduced coordinates and prediction from the SSM-reduced model. Subplot (\subref{fig:ross_rd}) shows the training and test data in three-dimensional reduced coordinates. Subplots (\subref{fig:ross_pre}) and (\subref{fig:ross_le}) show the time series forecasting results and the reconstructed Lyapunov exponent. The prediction remains accurate for about 3.3 Lyapunov times. The blue linear fit of the reduced-order model in the subplot (\subref{fig:ross_le}) gives the maximal Lyapunov exponent 0.0946, which is near 0.0938, the actual MLE of system (\ref{ross_eqa}).    }
    \label{fig:ross}
\end{figure}

As in the delay-embedded Lorenz attractor example in Section \ref{sec:3Dlorenz}, we observe only the $x$ coordinate of system (\ref{ross_eqa}) to generate scalar observable time series. Then, using the information of the fixed point $x_0$, we shift the origin of the system to this unstable fixed point. The training and test data are generated separately by the initial conditions $(1,1,1)$ and $(1,1,2)$ over the time interval $(0,7\,000)$ with a time step of $0.01$. After the first $(0,20)$ time interval, both the training and test trajectories lie very close to the Rössler attractor. Following the same procedure as in Section \ref{sec:3Dlorenz}, we use FNN to estimate the SSM dimension $d$ and pick the manifold expansion order according to the invariance error. We select the SSM attached to the unstable fixed point to be $d = 3$ dimensional, the delay-embedding dimension to be $\rho = 7$, and the polynomial expansion order to be $\mathcal{K} = 4$. Namely, we delay-embed the scalar observable in 7D and use \texttt{SSMLearn} to perform model reduction onto a 3D reduced space. The shape of the chaotic attractor in three-dimensional reduced coordinates is shown in Fig. \ref{fig:ross_rd}.

In the reduced coordinates, we use the kNN method to make forecasts for the test trajectory. We select $k = 4$ and perform the kNN prediction method on a training set containing $7 \times 10^{5}$ data points. From Fig. \ref{fig:ross}, the prediction stays accurate for around 3.3 Lyapunov times of system (\ref{ross_eqa}), and the reconstructed maximum Lyapunov exponent 0.0946 is close to the actual MLE of the system 0.0938. 


As in the other examples, we also tried to apply SINDy to the same seven-dimensional training data. We tested different combinations of SINDy parameters (polynomial order ranging from three to five, and the sparsity parameter $\lambda _{\text{SINDy}}$ ranging from $0.001$ to $0.1$), but failed to identify a feasible model. The prediction results generated by these SINDy models generally blew up, i.e., the models do not contain a stable attractor. This may be because the dynamics in these delay-embedded coordinates do not have a sparse polynomial form in general.

\subsection{Kuramoto–Sivashinsky equation} \label{sec:kse}

To test SSM-based model reduction in higher-dimensional chaotic systems, we consider the Kuramoto-Sivashinsky (KS) equation \cite{ChaosBook,kse}. This equation is often used as a simplified model to study spatiotemporal chaos and pattern formation in fluid dynamics, particularly in the context of understanding turbulence. Its dynamics in one spatial dimension is given by the PDE
\begin{equation} \label{kse_eqa}
    u_t + \frac{1}{2}(u^2)_x + u_{xx} + u_{xxxx} = 0, \;\;\; x\in \left[ -\frac{L}{2}, \frac{L}{2} \right], 
\end{equation}
where $x$ is the spatial coordinate and $t\geq 0$ is the time. Equation (\ref{kse_eqa}) describes the time evolution of the velocity $u = u(x,t)$ on a periodic domain $u(x,t) = u(x+L,t)$ with length $L$. When the velocity is $0$ everywhere on $x$, i.e., $u(x) \equiv 0$, the system has an unstable fixed point from which SSMs originate.

\begin{figure}[b]
    \centering
    \includegraphics[width=0.5\textwidth]{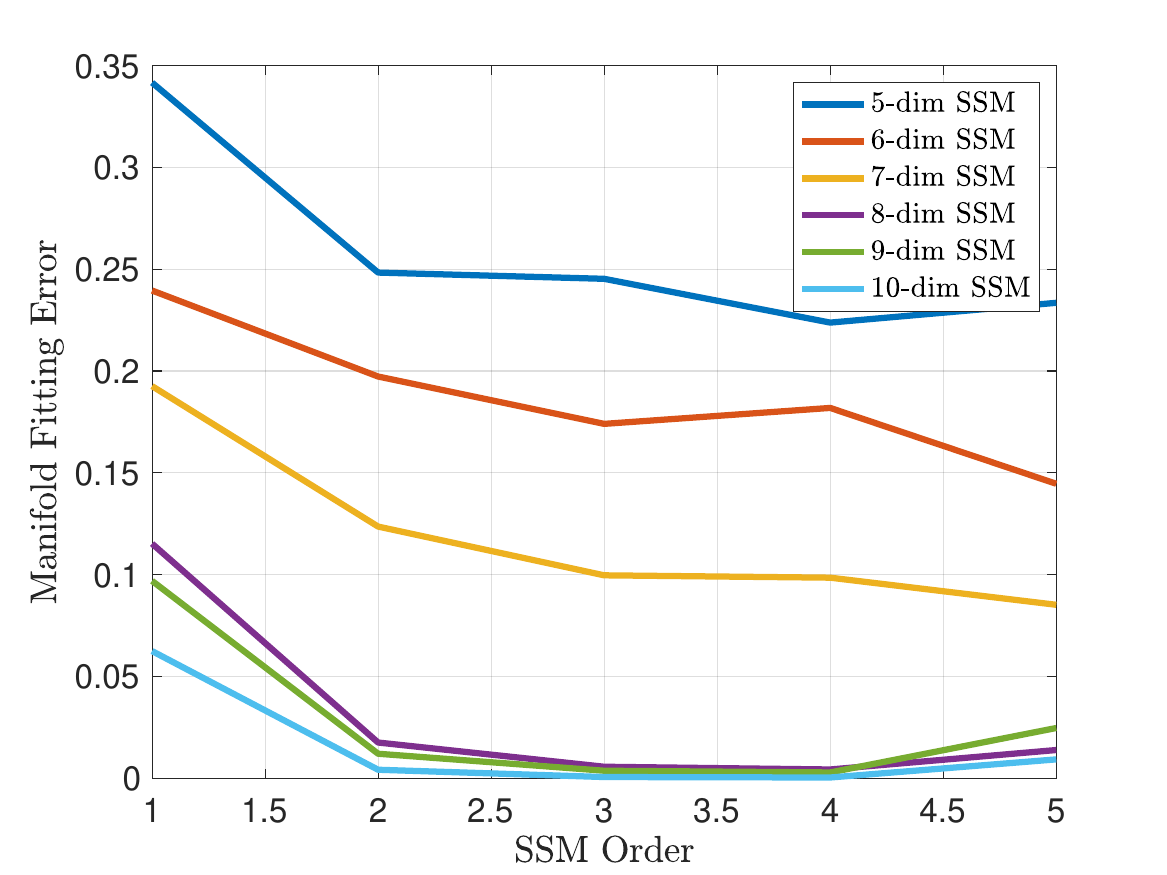}
    \caption{Manifold fitting error for different SSM dimensions, as functions of the SSM fitting order. From this figure, we fit a $d = 8$ dimensional SSM with a $\mathcal{K}=3$rd order expansion to the $1,024$-dimensional input data, which gives an invariance error $0.35 \%$ for the test set. }
    \label{fig:kse_error}
\end{figure}

To generate a data set for the discovery of SSMs, we adapt the \textsc{Matlab} code from \cite{ChaosBook} which uses a spectral solver for numerical solutions of the Kuramoto-Sivashinsky equation. We consider the system size $L=22$, for which a truncation of the number of Fourier modes to the range $16 \leq N \leq 128$ yields accurate solutions for this system size \cite{kse}, therefore we pick $N = 64$. We discretize $L$ into 1,024 segments, generate trajectories with the time step $0.25$, and discard the first 100 time units to get trajectories close to the attractor. The training set which is generated for the time length $1.25 \times 10^5$ contains approximately $5 \times 10^5$ data points. In this example, we only test the short-term predictability of the SSM-reduced model and therefore generate a short test set for a $100$ time interval.

For \texttt{SSMLearn}, computing the parametrization of an SSM requires solving the implicit minimization in equation (\ref{minimization}), which can be computationally demanding for high-dimensional systems. Therefore, we turn to the \texttt{fastSSM} algorithm \cite{Axs2022}, which significantly reduces the computational cost when the observable space dimension is high. The \texttt{fastSSM} algorithm calculates the tangent space of the SSM by singular value decomposition and the manifold parametrization via explicit polynomial regression. As a result, the computational cost scales linearly with the input data dimension. We refer to \cite{Axs2022} for a detailed discussion. 

\begin{figure}[t]
    \centering
    \includegraphics[width=0.7\textwidth]{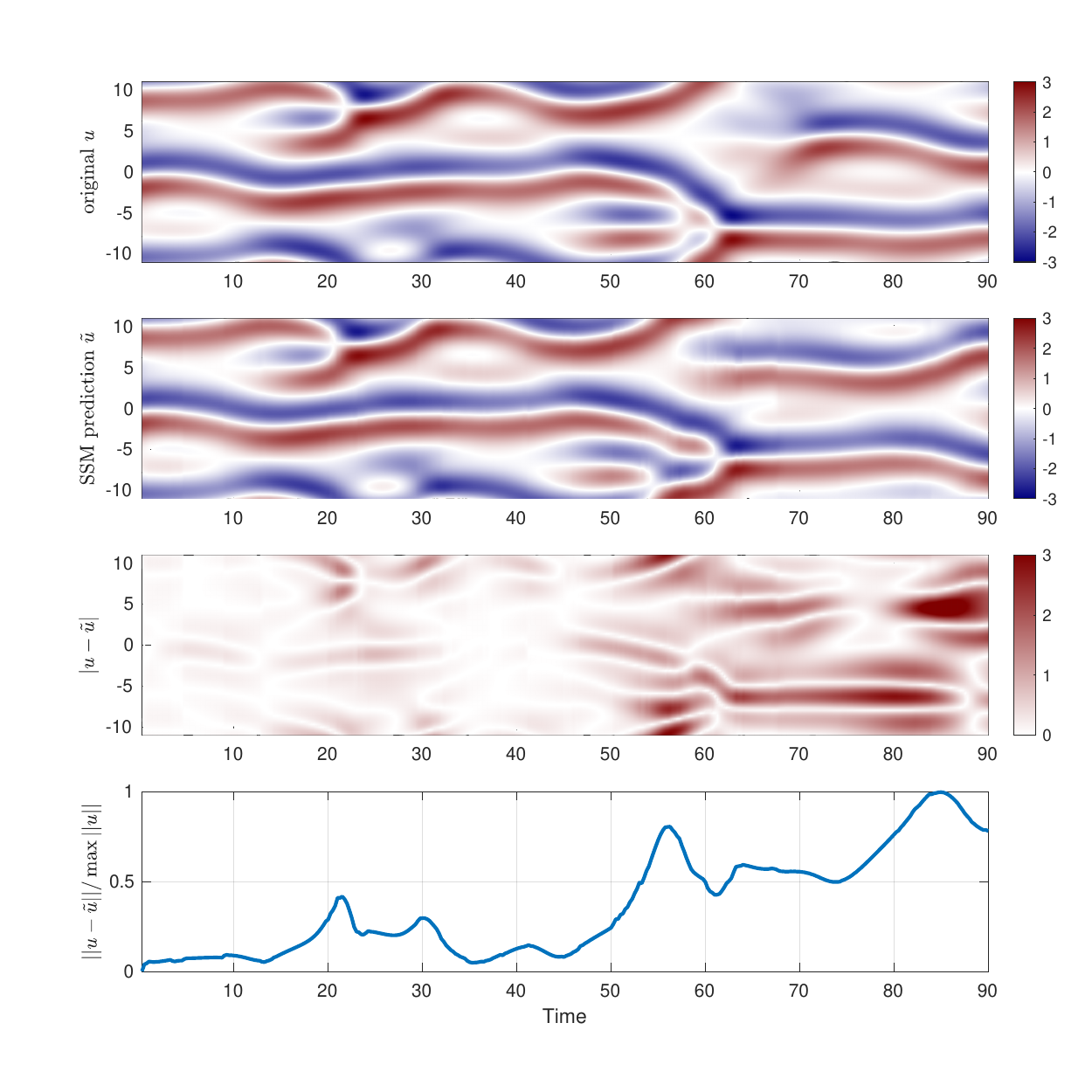}
    \caption{Test trajectory and its prediction by the SSM-reduced model with system size $L=22$. The top two plots are the velocity $u = u(x,t)$ and its prediction $\tilde{u} = \tilde{u}(x,t)$ for the time interval $(0,90)$. The SSM prediction resembles the real chaotic flow pattern but gradually diverges from it. The bottom two plots show the pointwise prediction error and the time evolution of the normalized trajectory error.}
    \label{fig:kse_pre}
\end{figure}

To estimate the SSM dimension from data, we plot the manifold fitting error as a function of the SSM dimension and SSM polynomial expansion order as shown in Fig. \ref{fig:kse_error}. When we increase the SSM dimension to eight, there is a significant drop in manifold fitting error. This indicates the chaotic attractor of the system can be captured by an inertial manifold with a dimensionality of at least eight. Therefore, we fit an eight-dimensional SSM with a $3$rd order polynomial expansion to the 1,024-dimensional trajectory data. The invariance error of the manifold fitting is $0.35 \%$ for the test data. The estimated dimension $d=8$ agrees with the estimations in \cite{PhysRevLett.117.024101, kse} and other model reduction results on inertial manifolds \cite{10.1063/5.0069536, Vlachas2022}.

After obtaining the data in 8D reduced coordinates, we perform the kNN algorithm with $k=9$ to make forecasts for the test trajectory. The forecasting quality and the validity range of each prediction depend on the amount of training data used and the individual initial condition of the test trajectory. We only illustrate one of the SSM prediction results in Fig. \ref{fig:kse_pre}. In general, the normalized trajectory errors reach 0.5 between 10 and 60 time intervals. A detailed analysis of Lyapunov exponents and average predictability for this infinite-dimensional example is beyond the scope of the present paper.


\subsection{Periodically forced finite-element beam} \label{sec:vkbeam}

\begin{figure}[h!]
    \centering
    \includegraphics[width=0.4\textwidth]{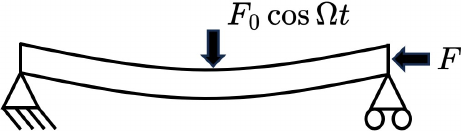}
    \caption{Dynamic buckling of a pinned-pinned von Kármán beam, with the subsequent addition of a time-periodic vertical load at its midpoint.}
    \label{fig:vk_beam}
\end{figure}

In our final and most complex example, we use data-driven model reduction to a mixed-mode SSM to identify a low-order model for the chaotic attractor of a periodically forced buckling beam. The nonlinear finite-element model we use for a von Kármán beam is a second-order approximation of the kinematics of a beam, taking into account the deformation of the cross-section for both bending and stretching. We specifically adopt here the finite-element \textsc{Matlab} code developed in \cite{vk} and the same buckling beam setting as in \cite{mixedmodeSSM}. The boundary condition of the beam is set as pinned-pinned, as shown in Fig. \ref{fig:vk_beam}, and an axial compressive force is applied to the right node.

The finite-element model of the von Kármán beam consists of 13 nodes, each with three degrees of freedom (DOFs): axial direction, transverse direction and rotation. Under the above boundary conditions, the beam model has 12 elements and 36 degrees of freedom. With the displacement and velocity of each DOF counted, the full model is a 72-dimensional dynamical system in its phase space. We adopt the physical parameters listed in Table \ref{tab:sys_param} from Ref. \cite{vk}.

\begin{table}[h]
\begin{center}
\begin{tabular}{cccc}\hline
Symbol & Meaning & value & [Unit]  \\ \hline
$t$ & time & - & [s] \\
$L$ & length of the beam & 2 & [m] \\
$h$ & height of the beam & 1e-2 & [m] \\
$b$ & width of the beam & 5e-2 & [m] \\  
$A$ & width of the beam & 5e-4 & [$\mathrm{m}^2$]  \\
$E$ & Young's modulus & 190e9 & [Pa]   \\  
$\kappa$ & viscous damping rate of material & 7e6 & [Pa$\cdot$s]   \\
$\rho$ & density & 7,850 & [kg/$\mathrm{m}^3$]   \\ \hline
\end{tabular}
\caption{Physical parameters for the von Kármán beam} \label{tab:sys_param}
\end{center}
\end{table}

Euler's critical load for buckling is given by $ P_{n}=\frac{n^{2}\pi^{2}EI}{l^{2}} $. For the first buckling mode $n=1$, the critical buckling force is
\begin{equation}
    F_{\mathrm{buckling}}=\frac{\pi^{2}EI}{l^{2}}=\frac{\pi^{2}Ebh^{3}}{12l^{2}}.
\end{equation}

For the buckling force value $F = 1.1F_{\mathrm{buckling}}$, the system has three fixed points: one unstable fixed point corresponding to purely axial displacement, and two stable fixed points corresponding to the upper and lower buckled states. The first 10 eigenvalues of the unstable fixed points are listed in Table \ref{vkbeam_eigenvalues}. Therefore, each unstable steady state has one pair of real eigenvalues $11.06$ and $-11.10$, the other pairs correspond to stable spirals. The positive real eigenvalue $11.06$ has the physical meaning of the beam leaving the unstable fixed point and approaching the upper or lower buckled states, and the negative real eigenvalue $-11.10$ corresponds to a trajectory settling down at the unstable state. The other complex pairs of eigenvalues represent the dynamics of higher modes. According to mixed-mode SSM theory \cite{mixedmodeSSM}, there exists a two-dimensional SSM tangent to the spectral subspaces corresponding to the two real eigenvalues. As in our earlier examples, these eigenvalues will not be used in our SSM construction.

\begin{table}[h]
\centering
\begin{tabular}{cccccc}
\toprule
$\lambda_{1}$ & $\lambda_2$ & $\lambda_{3,4}$ & $\lambda_{5,6}$ & $\lambda_{7,8}$ & $\lambda_{9,10}$ \\ \hline
11.06  & -11.10 & $-0.36\pm119.36i$ & $-1.83\pm295.56i$  & $-5.80\pm541.50i$ &   $-14.19\pm858.19i$  \\       \bottomrule               
\end{tabular}
\caption{First ten eigenvalues of the buckled von Kármán beam for reference.}
\label{vkbeam_eigenvalues}
\end{table}

\begin{figure}[h]
    \centering
    \includegraphics[width=0.5\textwidth]{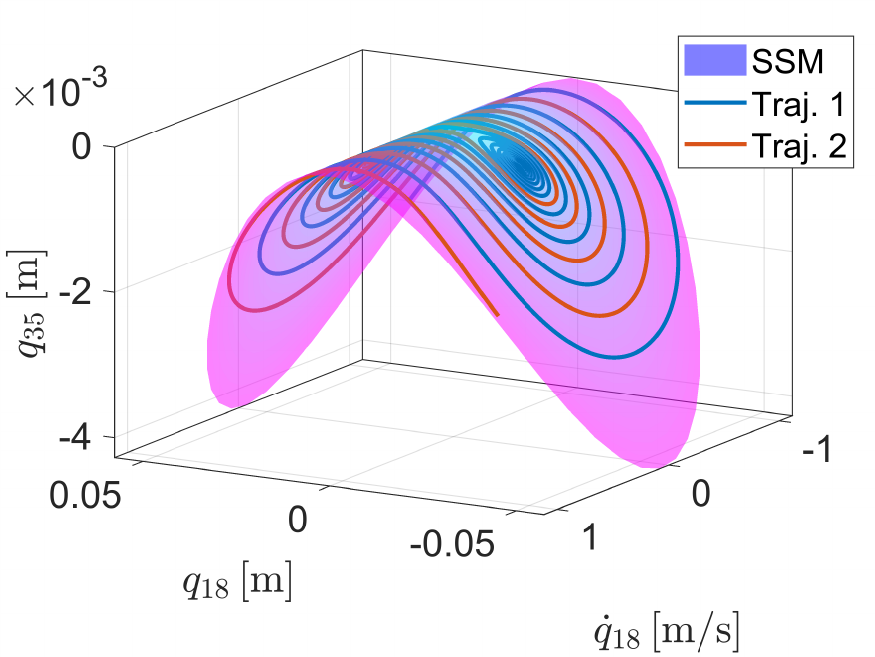}
    \caption{Non-chaotic training trajectories of the buckling beam model on a 2D, autonomous SSM.}
    \label{fig:vk_ssm}
\end{figure}

We generate four training trajectories from the autonomous beam model (without the vertical forcing shown in Fig. \ref{fig:vk_beam}) for the non-dimensionalized time interval $(0, 7.89)$, which covers $150$ times the oscillating period of the slowest mode. Each trajectory contains $1.5 \times 10 ^4$ points. The initial conditions are generated either by adding a small perturbation to the unstable stationary state, or by applying a large transverse force ($10^{4}$ N) at the middle node. After applying this force, the beam is further down than the lower buckled state, and will first oscillate around the two stable fixed points with a large amplitude, then converge to one of them. The higher modes corresponding to the complex eigenvalue pairs die out, and the trajectories converge to a 2D mixed-mode SSM, which we locate from the training trajectories using \texttt{SSMLearn}. Converged portions of two training trajectories and the SSM are shown in Fig. \ref{fig:vk_ssm}. Based on invariance error and NMTE for non-chaotic test data (one trajectory) in Table. \ref{vkbeam_order}, we construct a 7th-order SSM expansion and a 9th-order polynomial representation of the reduced dynamics in \texttt{SSMLearn}.

\begin{table}[]
\centering
\begin{tabular}{c|ccccc}
\hline
Manifold expansion order $(\mathcal{K})$                                                                    & 3      & 4      & 5      & 6      & 7      \\ \hline
Invariance error                                                                             & 0.10\%  & 0.10\% & 0.14\% & 0.14\% & 0.097\% \\ \hline
\begin{tabular}[c]{@{}c@{}}Reduced model order\\ (SSM expansion order $\mathcal{K}$ = 7)\end{tabular} & 7      & 8      & 9      & 10     & 11     \\ \hline
NMTE for non-chaotic test data                                                               & 0.24\% & 0.24\% & 0.12\% & 0.12\% & 0.16\% \\ \hline
\end{tabular}
\caption{Different choices of SSM expansion order and their corresponding invariance error for non-chaotic test data of the buckling beam. After selecting the manifold expansion order to be $\mathcal{K} = 7$, we choose the SSM-model polynomial order to be 9th according to their NMTE.}
\label{vkbeam_order}
\end{table}

Next, we apply transverse periodic forcing on the middle node of the beam, as shown in Fig. \ref{fig:vk_beam} which makes the beam oscillate chaotically. For moderate forcing amplitudes, the autonomous SSM, $\mathcal{W}(E)$, is known to persist in the form of a nearby, time-periodic SSM, $\mathcal{W}(E,t)$ \cite{mixedmodeSSM}. The leading-order forced dynamics on $\mathcal{W}(E,t)$ can be obtained by simply adding to the reduced dynamics of $\mathcal{W}(E)$ the projection of the time-periodic forcing on $E$ \cite{mixedmodeSSM}. This fact enables one to make predictions for the forced response based solely on a model trained on unforced data. 

\begin{figure}[h]
    \centering
    \begin{subfigure}{0.4\textwidth}
        \includegraphics[width=0.9\textwidth, height=2in]{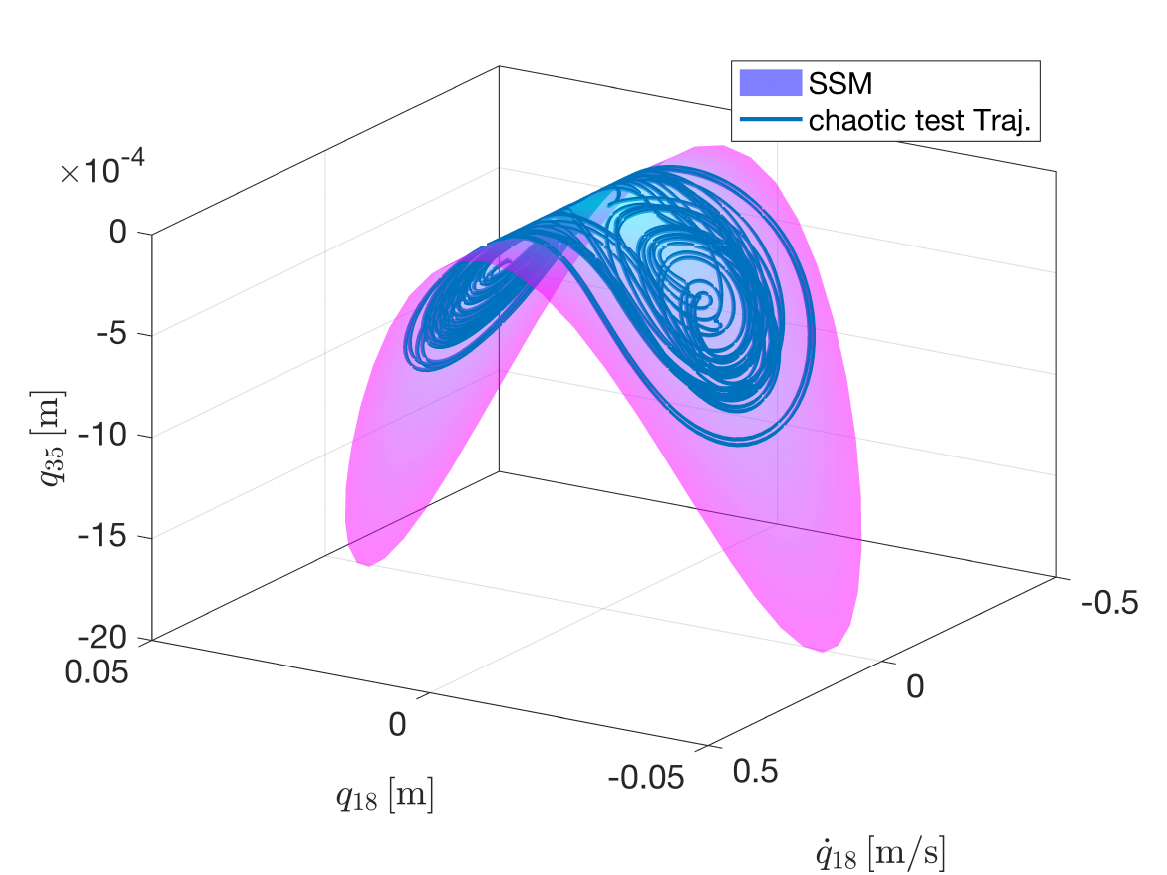}
        \caption{\label{fig:vk_ssm_chaotic}}
    \end{subfigure}
    \begin{subfigure}{0.4\textwidth}
        \includegraphics[width=0.9\textwidth, height=2in]{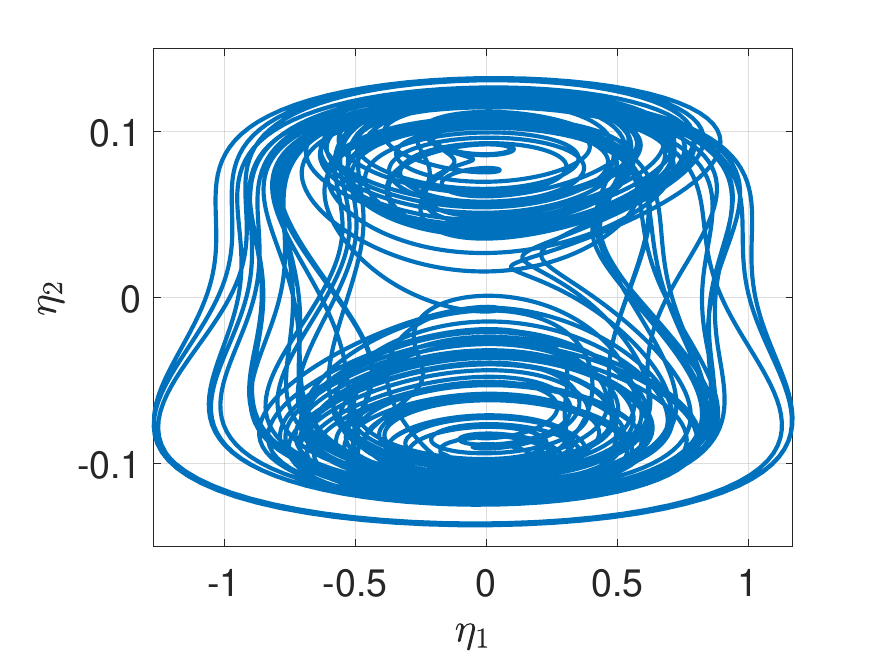}
        \caption{\label{fig:vk_rd}}
    \end{subfigure}
    \caption{(\subref{fig:vk_ssm_chaotic}) Chaotic response of the beam under periodic vertical forcing in the full phase space. This chaotic trajectory lies on the 2D SSM that we computed earlier from the non-chaotic data.
    (\subref{fig:vk_rd}) Chaotic trajectories of one testing data in two-dimensional reduced coordinates.}
    \label{fig:vk_ss}
\end{figure}

This procedure is implemented in \texttt{SSMLearn} \cite{SSMLearn,ssmlearn_code} and will be discussed here. The periodic forcing in Fig. \ref{fig:vk_beam} has the amplitude $F_0 = 3.4$ N and the frequency $\Omega = 15.4$ rad/s. We then generate the chaotic response of the buckling beam to this periodic force, and this one chaotic trajectory as the test data for our reduced-order modeling. Figure \ref{fig:vk_ss} shows this test trajectory lies entirely on the 2D SSM, and can be projected onto the two-dimensional reduced coordinates.

Specifically, up to leading order, the forced, SSM-reduced model on 2D SSM, $\mathcal{W}(E,t)$, is in the form (\cite{mixedmodeSSM})
\begin{equation} 
    \dot{\bm{\eta}} = \bm{R}(\bm{\eta}) + \begin{pmatrix}
    A \cos{\Omega t} \sin{\phi} \\
    A \cos{\Omega t} \cos{\phi} \end{pmatrix}.
\end{equation}

\begin{figure}[h!]
    \centering
    \begin{subfigure}{0.4\textwidth}
        \includegraphics[width=0.9\textwidth, height=2in]{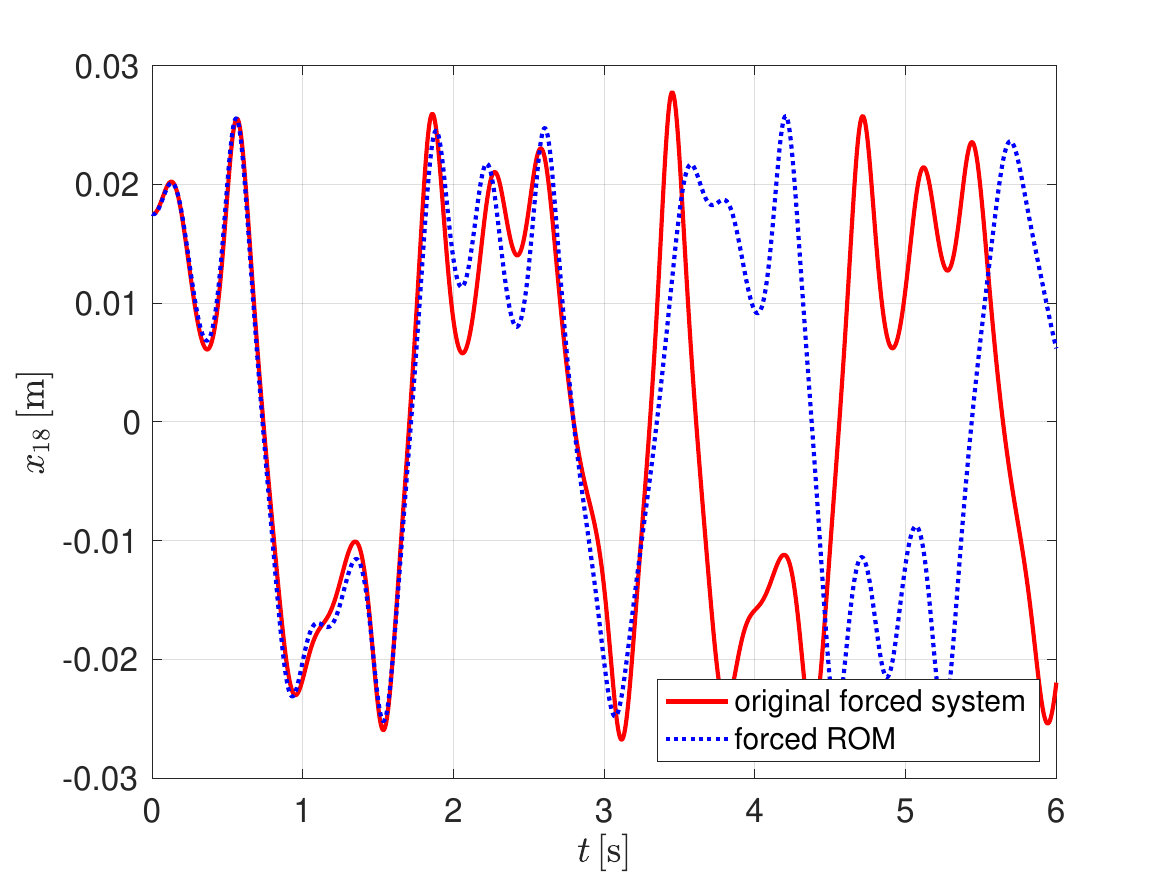}
        \caption{\label{fig:vk_pre}}
    \end{subfigure}
    \begin{subfigure}{0.4\textwidth}
        \includegraphics[width=0.9\textwidth, height=2in]{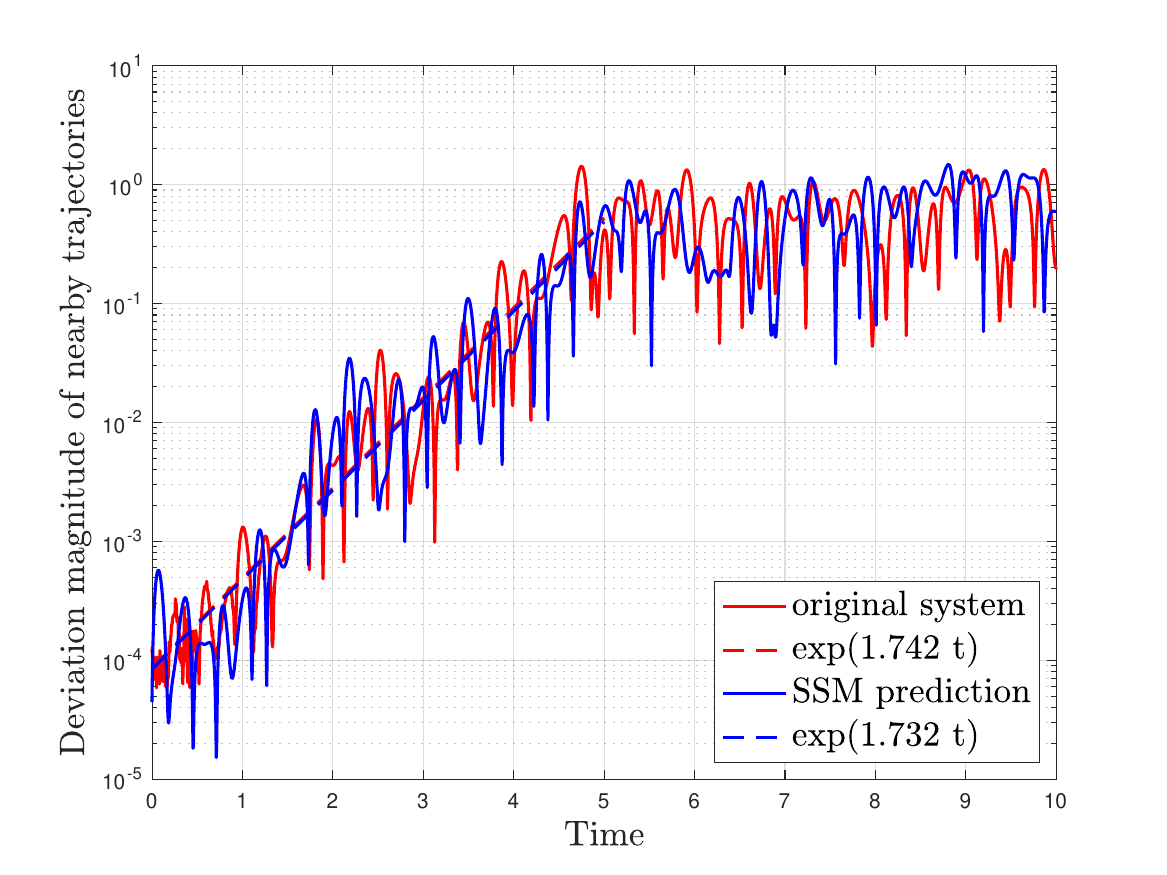}
        \caption{\label{fig:vk_ly}}
    \end{subfigure}
    \caption{(\subref{fig:vk_pre}) Prediction for a chaotic beam trajectory from an SSM-reduced model trained on unforced, non-chaotic data, then appended with the leading-order forcing term. The prediction lasts for about 3.31 Lyapunov time.
    (\subref{fig:vk_ly}) Separation of nearby trajectories against time. The original 72-dimensional system (red) and the reduced 2-dimensional model (blue) exhibit nearly equal Lyapunov exponents (1.742 and 1.732). } 
    \label{fig:vk_predict}
\end{figure}

After calibrating the forcing amplitude $A$ and forcing phase $\phi$ to a single forced trajectory (see \cite{SSMLearn} for details), we find that the SSM-reduced model computed from the autonomous system is indeed able to forecast the chaotic forced response, as shown in Fig. \ref{fig:vk_predict}. The Lyapunov exponent of the reconstructed system $1.732$ is very close to its counterpart $1.742$ calculated for the full, 72-dimensional system. Thus, the statistical properties of the chaotic attractor are successfully predicted by an SSM-reduced model in this example, too.

Even though a short-term prediction shown in Fig. \ref{fig:vk_pre} quickly loses its accuracy, one should not forget: This prediction by the 2D, SSM-reduced model is made entirely based on unforced, non-chaotic trajectories and yet it does approximate both the geometry and the long-term statistics (as described by the maximal Lyapunov exponent) with notable accuracy for the chaotic attractor of the full, 72-dimensional system. This shows the high physical interpretability of SSM-reduced models, which other data-driven model reduction techniques generally do not provide. For instance, the SINDy algorithm we used in some of our previous examples was not realistic to run on this high-dimensional example and the algorithm would not be equipped to make forced predictions from unforced training data.

\section{Discussion}

Using the recent theory of mixed-mode spectral submanifolds \cite{mixedmodeSSM}, we have developed a methodology to derive SSM-reduced models for dynamical systems with chaotic attractors. Specifically, we combined the SSMLearn algorithm with the method of nearest neighbors to obtain an accurate enough representation of the SSM-reduced dynamics on large enough domains that capture the attractors of the system. In order to contain the chaotic attractor and the dominant behavior of the underlying system, the SSM must have a dimension higher than the attractor dimension. In a strictly data-driven setting, the SSM dimension can be approximated using the false nearest neighbor method, then further determined by computing the invariance error.

This methodology is generally applicable to systems in which an inertial manifold containing the chaotic attractor is, in fact, an SSM emanating from an unstable steady state. In the examples treated here, that steady state was either a fixed point or a periodic orbit. Building the inertial manifold specifically as a spectral submanifold enables us to construct mathematically justified, low-dimensional models for the dynamics on and near the attractors. The models are either explicit polynomial ODEs or ODEs with right-hand sides interpolated using nearest neighbor points on the SSM. The core part of our approach is the open-source \texttt{SSMLearn} algorithm \cite{ssmlearn_code}, a general \textsc{Matlab} (Python) package for SSM-based model reduction. We additionally deploy the false nearest neighbor (FNN) algorithm \cite{fnnchaos} to identify the optimal SSM dimension from data. For SSM-based dimension reduction, no information about the time derivative data is required, and the amount of data needed is relatively low. 

We have applied two methods to model the reduced dynamics: a global polynomial fit and the kNN method. Both methods have limitations. The global polynomial fit requires computing the time derivatives of the trajectories and can become unstable for complex systems. The kNN method only uses local information to make predictions and therefore requires a large amount of data that densely covers the whole attractor. To counter these challenges, it is possible to combine SSM model reduction with other prediction methods, provide the weighting functions \cite{Crutchfield1987EquationsOM,PhysRevA.44.3496}, use other function bases such as radial basis functions \cite{CASDAGLI1989335}, or employ machine learning and deep learning methods \cite{KARUNASINGHE200692,Sangiorgio2022} such as neural ODEs \cite{10.1063/5.0069536}.

Our SSM-based attractor models have yielded accurate local predictions over short Lyapunov times and also reproduced closely the statistical and ergodic properties of chaotic attractors over longer times. Examples supporting these conclusions included the classic Lorenz and Rössler system, the extended 9D Lorenz model, a forced oscillator chain, the Kuramoto–Sivashinsky equation, and a periodically forced buckling beam. 

We also compare our approach with the data-driven nonlinear modeling method SINDy \cite{sindy}. While the advantages of SINDy include simplicity, ease of implementation and versatility, the method scales unfavorably with the input data dimension, and depends on the choice of the function library and the coordinates, as discussed in Section \ref{sec:rossler}. In efforts to address these challenges, SINDy has also been combined with various coordinate-identifying methods and model reduction techniques, such as learning appropriate observables in an assumed latent space. For SSM-based model reductions, the existence of such low-dimensional invariant manifolds is theoretically justified in the presence of a hyperbolic steady state. The requirement of such a steady state is a limitation of our method, but in our view, it is necessary for the localization of invariant manifolds from data to be on firm theoretical ground.

Importantly, since SSM-reduced models are fully dynamics-based, they can predict forced chaotic responses based solely on the knowledge of unforced (non-chaotic) trajectory data as demonstrated by our forced buckling beam example. We are unaware of any other data-driven model reduction approach that has this capacity.


Future applications of our methodology will include other PDE models, such as the 2D Kolmogorov flow \cite{kol1, kol2}, which also has finite-dimensional attractors and coexisting unstable steady states. Further extensions of mixed-mode SSM theory to more general anchoring steady states (such as invariant tori and other bounded invariant sets) are needed for the present methodology to become applicable to inertial manifolds that emanate from invariant sets other than fixed points or periodic orbits.

\section*{Data availability}
The code used to generate the numerical results included in this paper is available as part of the open-source \textsc{Matlab} package \texttt{SSMLearn}, available at \url{https://www.dropbox.com/scl/fi/8ydyh4vhxomkrx1o73im2/SSMLearn_chaotic.zip?rlkey=pjmz474e6x4ebkew08kuwx7jb&dl=0}.

\section*{Acknowledgments}
We acknowledge helpful discussions with Mattia Cenedese, Mingwu Li and Shobhit Jain.

\section*{Author Declearations}
The authors declare no conflict of interest.

\bibliography{aipsamp}

\end{document}